\documentclass[12pt,english]{article}
\usepackage[T1]{fontenc}
\usepackage[utf8]{inputenc}
\usepackage{babel}
\usepackage{geometry}
\usepackage{amsmath}
\usepackage{graphicx}
\usepackage{amssymb}

\makeatletter

\newtheorem{theorem}{Theorem}[section]
\newtheorem{corollary}[theorem]{Corollary}
\newtheorem{lemma}[theorem]{Lemma}

\numberwithin{equation}{section}

\makeatother

\newcommand{\eqnsection}{
\renewcommand{\theequation}{\thesection.\arabic{equation}}
 \makeatletter   \csname  @addtoreset\endcsname{equation}{section}
   \makeatother}

\def\proof{\noindent\textbf{Proof: $\,$} }
\def\qed{$\Box $}

\def\R{\mathbb{R}}

\def\Z{\mathbb{Z}}

\def\E{\mathbb{E}}
\def\P{\mathbb{P}}
\def\0{\mathbf{0}}
\def\1{\mathbf{1}}

\def\ol{\overline}

\renewcommand{\phi}{\varphi}

\begin{document}

\title{Spatial Brownian motion 
in renormalized Poisson potential: A critical case}

\date{March 3, 2011}

\author{Xia Chen%
\thanks{Supported in part by NSF grant DMS-0704024.%
} \, \, Jan Rosi\'nski%
\thanks{Supported in part by NSA grant MSPF-50G-049.%
}}

\maketitle

\def\theequation{$\ast$}

\begin{abstract}

Let $B_s$ be a three dimensional Brownian motion and
$\omega(dx)$ be an independent Poisson field on $\R^3$. It is proved that for any $t>0$, conditionally on $\omega(\cdot)$, 
\begin{equation} \label{*}
\E_0\exp\left\{\theta\int_0^t\ol{V}(B_s)ds\right\} \ \begin{cases}
< \infty \ a.s.   & \text{if }  \theta< 1/16, \medskip \\ 
 =  \infty \ a.s.   & \text{if }  \theta> 1/16, 
\end{cases}
\end{equation}
where $\ol{V}(x)$ is the renormalized Poisson potential
$$
\ol{V}(x)=\int_{\R^3}{1\over\vert x-y\vert^2}\big[\omega(dy)-dy\big].
$$
Then the long term behavior of the quenched exponential moment  \eqref{*} is determined for  $\theta \in (0, 1/16)$ in the form of integral tests.

This paper exhibits and builds upon the interrelation between the exponential moment  \eqref{*}  and the celebrated Hardy's inequality
$$
\int_{\R^3}{f^2(x)\over\vert x\vert^2}dx\le 4\|\nabla f\|_2^2, 
\hskip.2in f\in W^{1,2}(\R^3).
$$
\end{abstract}

\begin{quote} {\footnotesize

}

\end{quote}

\newpage

\renewcommand{\theequation}{\thesection.\arabic{equation}}

\section{Introduction } \label{intro}

Consider a particle moving randomly according to a standard $d$-dimensional Brownian motion $B_s$ in $\R^d$.   Independently, there is a family of obstacles
randomly placed in the space $\R^d$ according to a Poisson field $\omega(dx)$ (i.e., a Poisson random measure). Assume that each obstacle has mass 1 and the Poisson field $\omega(dx)$ has the Lebesgue measure $dx$ as its intensity measure. Throughout this paper, ``$\P_z$'' and ``$\E_z$'' will respectively stand for the probability law and expectation relative to Brownian motion $B_s$ with $B_0=z$. Notation ``$\P$'' and ``$\E$'' will be used for the 
probability law and expectation, respectively, 
relative to  Poisson field $\omega(dx)$. 

Given a {\em shape function} $K(x)$ (as known by mathematicians, and a {\em point-mass potential} by physicists) on $\R^d$, the potential associated with the random 
mass distribution $\omega(\cdot)$
is given by 
$$
V(x) = \int_{\R^d} K(x-y) \, \omega(dy),
$$
and is called a {\em Poisson potential}.  
The quantity
$$
t^{-1} \int_0^t V(B_s) \, ds
$$
represents the average Poisson potential along a Brownian trajectory. 
More important quantities of interest are
the respective {\em annealed} and {\em quenched} exponential moments
\begin{align}\label{intro-0}
\E\otimes \E_0\exp\bigg\{\pm \int_0^t V(B_s) \, ds\bigg\}\hskip.1in\hbox{and}
\hskip.1in
\E_0\exp\bigg\{\pm  \int_0^t V(B_s) \, ds\bigg\}.
\end{align}
Knowledge of their asymptotic behavior at $t \to \infty$ is fundamental to our understanding of parabolic Anderson models (see Corollary \ref{anderson}).
The reader is referred to  \cite{BK2}, \cite{CM}, \cite{CV},
\cite{CGM}, \cite{DM}, \cite{DV}, \cite{FV}, \cite{GK}, \cite{GKM}, \cite{G-M}, \cite{Pastur77}, \cite{Povel}, \cite{Stolz}, and \cite{Sznitman} for the existing literature on this topic.

In the classical literature on this subject, the function $K(x)$ of the Poisson potential was assumed to be bounded and/or compactly supported.
However, in physics many point-mass potential functions are unbounded.
For example, in scattering theory, power potentials  
$K(x)= \mp |x|^{-p}$ ($d=3$) play  a significant role, 
see \cite{McDaniel},   \cite{Newton}. The parameter $p>0$ is called the 
index of attraction or repulsion, respectively. When $p=1$, we have Coulomb interaction. $K(x)=|x|^{-4}$ is referred to as Maxwellian potential, $K(x)=-|x|^{-4}$ is important in the study of ionized gases, $K(x)=|x|^{-2}$ is known as a centrifugal potential, see \cite[Chapters 1-7, 3-6]{McDaniel}.

Donsker and Varadhan \cite{DV}, Pastur \cite{Pastur77}, and Fukushima \cite{Fukushima} 
studied the asymptotics in \eqref{intro-0} with the {\em negative} signs for the case $K(x)= |x|^{-p}$. Specifically, in \cite{DV} and \cite{Pastur77} the asymptotics of the annealed moment were obtained when $p>d+2$
and $d<p<d+2$, respectively. In \cite{Fukushima} both  annealed and quenched  moments are determined for $d<p<d+2$. In these papers the singularity of $K(x)$ was circumvented by applying truncations near the origin.

When $p\le d$, the Poisson potential becomes infinite a.s. 
To deal with this problem, in the recent paper \cite{C-K}, the renormalized Poisson potential
$$
\ol{V}(x)=\int_{\R^d}{1\over\vert x-y\vert^p} \, \big[\omega(dy)-dy\big]
\hskip.2in x\in\R^d
$$
was introduced; it exists as a random integral if and only if $d/2<p<d$, see Corollary 1.3 and
 physical arguments behind the renormalization in \cite{C-K}. In the same work, 
integrabilities associated with the construction
of the annealed and quenched exponential moments
$$
\E\otimes\E_0\exp\bigg\{\pm\theta\int_0^t\ol{V}(B_s) \, ds\bigg\}
\hskip.1in\hbox{and}
\hskip.1in
\E_0\exp\bigg\{\pm\theta\int_0^t \ol{V}(B_s) \, ds\bigg\}
$$
were investigated.  In the range $d/2<p<d$,  the time-integral
$$
\int_0^t\ol{V}(B_s) \, ds
$$
is well defined and satisfies the annealed integrability (and therefore
quenched integrability as well)
$$
\E\otimes\E_0\exp\bigg\{-\theta\int_0^t\ol{V}(B_s) \, ds\bigg\}<\infty
$$
for every $\theta>0$ and $t>0$. We refer the recent papers \cite{Chen-1}
and \cite{C-K-1} for the study on
the asymptotics of quenched and annealed negative exponential moments, respectively.

However, the case of exponential moments with positive
coefficient is far more delicate.
By  \cite[Theorem 1.4]{C-K}, for every $\theta>0$ and $t>0$
$$
\E\otimes\E_0\exp\bigg\{\theta\int_0^t\ol{V}(B_s) \, ds\bigg\}=\infty.
$$
On the other hand, by \cite[Theorem 1.5]{C-K} the quenched exponential moment exists
for any $\theta>0$, $t>0$, and $p<2$, as we have with probability 1,
\begin{equation}\label{intro-1}
\E_0\exp\bigg\{\theta\int_0^t\ol{V}(B_s)ds\bigg\} \ \begin{cases}
< \infty \ a.s.  & \text{if }  p< 2, \medskip \\ 
= \infty \ a.s.   & \text{if }  p> 2. 
\end{cases}
\end{equation}
Furthermore, the first author recently observed in \cite{Chen-1} that
\begin{align}\label{intro-2}
\lim_{t\to\infty}  {1\over t} \Big({\log\log t\over\log t}\Big)^{2\over 2-p} &
\log\E_0\exp\bigg\{\theta\int_0^t\ol{V}(B_s)ds\bigg\}\\
&={1\over 2}p^{p\over 2-p}
(2-p)^{4-p\over 2-p}\Big({d\theta\sigma(d,p)\over 2+d-p}\Big)^{2\over 2-p}
\hskip.2in a.s. -\P, \nonumber
\end{align}
where $\sigma(d,p)>0$ is the best constant of the inequality
\begin{align}\label{intro-3}
\int_{\R^d}{f^2(x)\over\vert x\vert^p}dx\le C\|f\|_2^{2-p}\|\nabla f\|_2^p, 
\hskip.2in f\in W^{1,2}(\R^d).
\end{align}
The only unanswered case is $p=2$ and necessarily $d=3$ (recall the constraint ${d/2<p<d}$); we will call it the {\em critical case}.
Our results will justify this name.

The present paper is devoted to the study of the quenched exponential moment
\begin{equation} \label{intro-3.5}
\E_0\exp\bigg\{\theta\int_0^t\ol{V}(B_s) \, ds\bigg\}, \hskip.2in\theta>0,\hskip.1in
t>0
\end{equation}
in the critical case, i.e., when
\begin{equation} \label{intro-3.6}
\ol{V}(x)=\int_{\R^3}{1\over\vert x-y\vert^2} \, \big[\omega(dy)-dy\big], 
\hskip.2in x\in\R^3.
\end{equation}
Remark that in physics $K(x)= \pm |x|^{-2}$ ($x \in \R^3)$ is a transition potential. 
It lies on the boundary between the classes of regular ($p<2$) and singular ($p>2$) potentials separating 
fundamentally different physical systems, see  \cite[Section II]{FLS}. For example, in nonrelativistic quantum mechanics a particle in an attractive singular potential has infinite negative energy. The particle in this case "falls" to the center with infinite velocity. However, if $p<2$, the energy is finite, solutions to physical problems are uniquelly given,  and there is no problem with their physical interpretation, see  \cite[Sections I--II.A]{FLS}.

It was already noticed in \cite[Theorem 1.5]{C-K} that the quenched exponential moment in \eqref{intro-3.5}
is infinite a.s. for  $\theta$  sufficiently large and all $t>0$. A natural question is whether
this is true for all $\theta>0$. Fortunately, the answer is negative. If so, what is the critical value  
$\theta_0$ where the phase transition occurs? (It is even not clear that $\theta_0$ must be deterministic.)
We prove that $\theta_0=1/16$. Then we establish the asymptotic behavior of the quenched exponential moment in 
\eqref{intro-3.5}, showing that it is fundamentally different from \eqref{intro-2} since the strong law of large numbers does not hold in the critical case. These results are the consequences of the interrelation with Hardy's inequality \eqref{H0} via a chain of asymptotic equivalences sketched in \eqref{intro-13}.
In conclusion, the critical case of $p=2$ is substantially different from the other cases. The only continuity 
appears in Hardy's inequality, where a formal substitution of $p=2$ in \eqref{intro-3} gives \eqref{H0}.

The paper is organized as follows. In section \ref{main} we present main results and their application
to the parabolic Anderson model.  In section \ref{estimate},
we develop key tools for the estimations needed in later sections.
Some of these tools are interesting for their own novelty. 
Slepian-type correlation inequalities 
for infinite divisible fields (cf. \cite{S-T})  are 
provided (Lemma \ref{lem:00}) for the proof of (\ref{intro-15}), where the random
variables $\omega(z+Q_{b\delta})$ ($z\in 2\delta\Z^3$) are correlated
as $b>1$. An estimation by a chaining maximal inequality (Lemma \ref{poisson-2}) allows
the truncation of the Poisson potential at the proper level. Feynman-Kac
formula plays a crucial role in the proof of the main results in this paper.
A clean and simple
minorization bound (Lemma \ref{exit}) for Brownian density killed upon
exit
leads to a Feynman-Kac lower bound (Lemma \ref{FK-6}) adoptable to our setting.
For the Feynman-Kac upper bound (Lemma \ref{FK-1}) with the random potential
$\ol{V}(\cdot)$, we  use  the independence between the Brownian exit time
and Brownian exit location from a ball centered at 0.
The lower and upper bounds for the main theorems are proved in the sections
\ref{lb} and \ref{ub}, respectively. The main ingredients in these two sections
are the estimation of the principal eigenvalues of the correspondent
initial-boundary value problems that leads to
the relation suggested by (\ref{intro-13}) and, the strong laws for 
extreme values of the Poisson field indicated by (\ref{intro-14})--(\ref{intro-15}). 
Section \ref{H} is devoted to Hardy's inequality and related facts.

\bigskip

\section{Main results} \label{main}

From now on we will assume $d=3$, $p=2$ and that the renormalized Poisson potential $\ol{V}(x)$ is given by \eqref{intro-3.6}, if not otherwise stated. 

\begin{theorem}\label{intro-4} 
For every $t>0$,
\begin{equation}\label{intro-5}
\E_0\exp\bigg\{\theta\int_0^t\ol{V}(B_s)ds\bigg\} \ \begin{cases}
< \infty \ a.s.  & \text{if }  \theta< 1/16, \medskip \\ 
= \infty \ a.s.    & \text{if } \theta > 1/16. 
\end{cases}
\end{equation}

\end{theorem}

In view of the limit law (\ref{intro-2}) obtained in the non-critical
case, a natural problem is the asymptotic behaviors in the critical
case. Recall that a positive function $\gamma(t)$ on $\R^+$ is said to
be regularly varying at infinity if the limit
$$
\lim_{t\to\infty}{\gamma(\lambda t)\over \gamma(t)}=c(\lambda)
$$
exists for each $\lambda>0$. A regularly varying function
$\gamma(t)$ is said to be slowly varying at infinity, if
$c(\lambda)\equiv 1$. From Karamata theory, every regularly varying 
function $\gamma(t)$ has a representation
$\gamma(t)=t^\beta l(t)$, where $\beta$ is a constant and $l(t)$ is
a slowly varying function.

Throughout, $l(t)$ will stand for a slowly varying function at infinity.

\begin{theorem}\label{intro-9}
For every $\theta \in (0, 1/16)$
\begin{align}\label{intro-10}
\limsup_{t\to\infty} \ & t^{-{k+1\over k-1}}l(t)^{-{2\over 3(k-1)}} 
\log\E_0\exp\bigg\{\theta \int_0^t\ol{V}(B_s)ds\bigg\}  \\
& = \begin{cases}
 0 \ a.s.   & \text{if} \  \displaystyle\int_1^\infty {dt\over t \cdot l(t)}<\infty, \medskip \\ 
 \infty \ a.s.   & \text{if} \  \displaystyle\int_1^\infty {dt\over t \cdot l(t)}<\infty, 
\end{cases}   \nonumber 
\end{align}	
where $k=\lfloor(8\theta)^{-1}\rfloor$
 is the integer part of $(8\theta)^{-1}$.
\end{theorem}

\begin{theorem}\label{intro-11}
For every $\theta \in (0, 1/16)$  
\begin{align}\label{intro-12}
\liminf_{t\to\infty} \ & t^{-{k+1\over k-1}}l(t)^{-{2\over 3(k-1)}} 
\log\E_0\exp\bigg\{\theta \int_0^t\ol{V}(B_s)ds\bigg\}  \\
& = \begin{cases}
 0 \ \ a.s.   & \text{if} \  \displaystyle\int_1^\infty {1\over t }\exp\big\{-c\cdot l(t)\big\}dt
=\infty \ \text{for some $c>0$}, \medskip \\ 
 \infty \ a.s.   & \text{if} \  \displaystyle\int_1^\infty {1\over t }\exp\big\{-c\cdot l(t)\big\}dt
<\infty \ \text{for every $c>0$}, 
\end{cases}   \nonumber 
\end{align}
where $k=\lfloor(8\theta)^{-1}\rfloor$ is as in Theorem \ref{intro-9}.	
\end{theorem}

Theorems \ref{intro-9}--\ref{intro-11} show rather unexpected behavior of the quenched exponential moments with regard to $\theta$. Indeed, putting $\theta$ into different sub-intervals of the
partition
$$
\Big(0, {1\over 16}\Big)=\Big({1\over 24}, {1\over 16}\Big)\cup
\bigcup_{k=3}^\infty\Big({1\over 8(k+1)}, {1\over 8k}\Big]
$$
leads to different asymptotic rates. On the other hand, moving $\theta$ around within the same sub-interval 
does not bring any  change to the asymptotic behavior of the system.
\medskip

Our main results indicates that as far as the strong limit is concerned,
there is not ``right'' deterministic normalization to the logarithm of the
quenched exponential moment in the critical setting $p=2$ and $d=3$. Indeed,
by Theorems \ref{intro-9}--\ref{intro-11}, for any 
$\displaystyle 0<\theta< 1/16$, and for any
positive deterministic function $\gamma(t)$ regularly varying at infinity,
with probability 1
$$
\limsup_{t\to\infty}\gamma(t)^{-1}
\log\E_0\exp\bigg\{\theta\int_0^t\ol{V}(B_s)ds\bigg\}=0
\hskip.05in\hbox{or}\hskip.05in \infty
$$
and
$$
\liminf_{t\to\infty}\gamma(t)^{-1}
\log\E_0\exp\bigg\{\theta\int_0^t\ol{V}(B_s)ds\bigg\}=0
\hskip.05in\hbox{or}\hskip.05in \infty.
$$
This pattern sharply contrasts \eqref{intro-2} observed in the non-critical
setting.

Letting $l(t)$ be some specific functions, we get the following results:
$$
\limsup_{n\to\infty}t^{-{k+1\over k-1}}(\log t)^{-{2\over 3(k-1)}}
\log\E_0\exp\bigg\{\theta\int_0^t\ol{V}(B_s)ds\bigg\}=\infty\hskip.2in a.s.
$$
On the other hand, for any $\delta>0$
$$
\limsup_{n\to\infty}t^{-{k+1\over k-1}}\big((\log t)(\log\log)^{1+\delta}
\big)^{-{2\over 3(k-1)}}
\log\E_0\exp\bigg\{\theta\int_0^t\ol{V}(B_s)ds\bigg\}=0\hskip.2in a.s.
$$
As for the liminf behavior,
$$
\liminf_{n\to\infty}t^{-{k+1\over k-1}}(\log\log t)^{{2\over 3(k-1)}}
\log\E_0\exp\bigg\{\theta\int_0^t\ol{V}(B_s)ds\bigg\}=0\hskip.2in a.s.
$$
On the other hand, for any $l(t)\gg \log\log t$ as $t\to\infty$,
$$
\liminf_{n\to\infty}t^{-{k+1\over k-1}}l(t)^{{2\over 3(k-1)}}
\log\E_0\exp\bigg\{\theta\int_0^t\ol{V}(B_s)ds\bigg\}=\infty\hskip.2in a.s.
$$

Theorem \ref{intro-4} provides solution to the parabolic Anderson equation
\begin{align}\label{intro-6}
\left\{\begin{array}{ll}\partial_tu(t,x)
=\kappa\Delta u(t,x)+\theta \ol{V}(x)u(t,x)\\\\
u(0, x)=1\end{array}\right.
\end{align}
where $\kappa>0$ is a constant called diffusion coefficient.
Indeed, consider the time-space field
\begin{align}\label{intro-7}
u_\theta(t,x)=\E_x\exp\bigg\{\theta\int_0^t\ol{V}(B_{2\kappa s})ds\bigg\}
=\E_x\exp\bigg\{{\theta\over 2\kappa}\int_0^{2\kappa t}\ol{V}(B_{s})ds\bigg\}.
\end{align}
By translation invariance of
the Poisson field, for any $x\in\R^d$
\begin{align}\label{intro-8}
\Big\{u_\theta(t,x);\hskip.1in t\ge 0\Big\}\buildrel d\over =
\Big\{u_\theta(t,0);\hskip.1in t\ge 0\Big\}.
\end{align}
By Theorem \ref{intro-4}, $u_\theta(t,x)<\infty$ a.s. for every $x\in\R^d$
and $t>0$ when $\theta<\kappa /8$. The argument same as the one for 
Proposition 1.6, \cite{C-K} concludes that when $\theta<\kappa /8$,
$u_\theta(t,x)$ is a mild solution
to the equation (\ref{intro-6}) in the sense that
$$
\int_0^tp_{2\kappa (t-s)}(x-y)\vert\ol{V}(y)\vert u_\theta(s, y)dyds<\infty
\hskip.2in x\in\R^3,\hskip.1in t>0
$$
and
$$
u_\theta(t,x)=1+\theta\int_0^tp_{2\kappa (t-s)}(x-y)\ol{V}(y)u_\theta(s, y)dyds
\hskip.2in x\in\R^3,\hskip.1in t>0
$$
where $p_t(x)$ is the Brownian density.

Further, Theorem
\ref{intro-9} and Theorem \ref{intro-11} lead to the long term property
of the stochastic partial different equation \eqref{intro-6}.

\begin{corollary}\label{anderson} 
Under $d=3$ and $p=2$, the random field
$u_\theta(t,x)<\infty$ for all $\theta<\kappa/8$ and 
$(t,x)\in\R^+\times\R^3$, and 
$u_\theta(t,x)=\infty$ for all $\theta>\kappa/8$ and
$(t,x)\in\R^+\times\R^3$. When $\theta<\kappa/8$,
$u_\theta(t,x)$ is a mild solution to the equation (\ref{intro-6})
and further, for any $x\in\R^3$,
\begin{align}\label{anderson-1}
&\limsup_{t\to\infty}t^{-{i+1\over i-1}}l(t)^{-{2\over 3(i-1)}}
\log u_\theta(t,x)=\left\{\begin{array}{ll} 0\hskip.2in a.s.\hskip.05in\hbox{if}
\hskip.05in\displaystyle\int_1^\infty {dt\over t \cdot l(t)}<\infty\\\\
\infty\hskip.2in a.s.\hskip.05in \hbox{if}
\hskip.05in\displaystyle\int_1^\infty {dt\over t \cdot l(t)}=\infty
\end{array}\right.
\end{align}
\begin{align}\label{anderson-2}
\liminf_{t\to\infty} & t^{-{i+1\over i-1}}l(t)^{{2\over 3(i-1)}}
\log u_\theta(t,x)\\
&=\left\{\begin{array}{ll} 0\hskip.2in a.s.\hskip.05in\hbox{if}
\hskip.05in\displaystyle\int_1^\infty {1\over t }\exp\big\{-c\cdot l(t)\big\}dt
=\infty \hskip.05in\hbox{for some $c>0$}\\\\
\infty\hskip.2in a.s.\hskip.05in \hbox{if}
\hskip.05in\displaystyle\int_1^\infty {1\over t }\exp\big\{-c\cdot l(t)\big\}dt
<\infty \hskip.05in\hbox{for every $c>0$}
\end{array}\right.\nonumber
\end{align}
where $i=\lfloor (4\theta)^{-1}\kappa\rfloor$
is the integer part of $(4\theta)^{-1}\kappa$.
\end{corollary}

Given the non-deterministic asymptotic behaviors observed from
Theorem \ref{intro-9} and Theorem \ref{intro-11}, the weak law
(if any) becomes an interesting problem. In view of
Theorem \ref{intro-9} and Theorem \ref{intro-11}, one might expect
that the process
$$
t^{-{k+1\over k-1}}
\log\E_0\exp\bigg\{\theta\int_0^t\ol{V}(B_s)ds\bigg\}
$$
converges to a non-degenerated distribution. We leave this problem
to future study.

The critical ($p=2$) and non-critical ($p<2$) 
settings depend on  the environment in different ways and therefore
are treated differently.
In the non-critical case, the quantity
$$
\E_0\exp\bigg\{\theta\int_0^t\ol{V}(B_s)ds\bigg\}
$$
is made by letting Brownian particle stay in a slowly
shrinking neighborhood
that provides maximal energy from Poisson field among all same-size
neighborhoods in a large ball of the radius (roughly) $t$.
Consequently, the limit in (\ref{intro-2}) depends on the extreme
values of the Poisson potential $\ol{V}(\cdot)$ over a group of 
shrinking neighborhoods.

In contrary, the limsup in Theorem \ref{intro-9} and the liminf in
Theorem \ref{intro-11} correspond with the value $\infty$
to the existence (with a proper asymptotic intensity)
of the neighborhoods in which the number of Poisson 
obstacles exceeds the fixed level $k=[(8\theta)^{-1}]$ within proper
distance, and
with the value 0 to the absence of such neighborhoods.
The central piece behind this strategy is
the celebrated Hardy's inequality (Lemma \ref{H-0'}) which states that
\begin{equation} \label{H0}
\int_{\R^3}{f^2(x)\over\vert x\vert^2}dx\le 4\|\nabla f\|_2^2
\hskip.2in f\in W^{1,2}(\R^3).
\end{equation}
As a consequence (Lemma \ref{H-2})
of Hardy's inequality,
$$
H(\theta)\equiv\sup_{g\in {\cal F}_3}\bigg\{\theta\int_{\R^3}
{g^2(x)\over\vert x\vert^2}dx
-{1\over 2}\int_{\R^3}\vert\nabla g(x)\vert^2dx\bigg\}
=\left\{\begin{array}{ll} 0\hskip.2in\hbox{if}\hskip.05in\theta\le 1/8 \\\\
\infty\hskip.2in \hbox{if} \hskip.05in\theta> 1/8.
\end{array}\right.
$$
The connection of Theorem \ref{intro-9} and Theorem \ref{intro-11} to Hardy's inequality is
described roughly by the following almost sure
asymptotic relation:
\begin{align}\label{intro-13}
\log \E_0 & \exp\bigg\{\theta\int_0^t\ol{V}(B_s)ds\bigg\}\\
&\approx
\log \E_0\Bigg[\exp\bigg\{\theta\int_0^t\ol{V}(B_s)ds\bigg\};\hskip.05in
\max_{0\le s\le t}\vert B_s\vert\le R\Bigg]\nonumber\\
&\approx t^{k+1\over k-1}l(t)^{\pm{2\over 3(k-1)}}\bigg\{o(1)+
A(t)H\Big(\theta\max_{\stackrel{\scriptstyle\vert z\vert\le R}
{\scriptstyle z\in 2\delta\Z^3}}\omega (z+Q_{b\delta})\Big)\bigg\}\nonumber
\end{align}
where $R$ increases to $\infty$ and $\delta>0$ decreases to zero with
suitable polynomial rates as $t\to\infty$, where $z+Q_{b\delta}$ represents
the cubic $z+[-b\delta, b\delta]^3$ with $b$ being a fixed constant,
and where $A(t)$ ranges from a constant (in the argument for the lower bound)
to a function increasing to $\infty$ at a considerable speed (in the argument
for the upper bound).

Our strategy is to let the Brownian particle spend  significant portion of 
the duration $[0, t]$
in one of the $\delta$-neighborhoods within the distance $R$.
A  principle of choosing $R$ and $\delta$
is to make  alternation between the behaviors
\begin{align}\label{intro-14}
\limsup_{t\to\infty}\max_{\stackrel{\scriptstyle\vert z\vert\le R}
{\scriptstyle z\in 2\delta\Z^3}}
\omega (z+Q_{b\delta})\le (8\theta)^{-1}\hskip.2in a.s.
\end{align}
\begin{align}\label{intro-15}
\liminf_{t\to\infty}\max_{\stackrel{\scriptstyle\vert z\vert\le R}
{\scriptstyle z\in 2\delta\Z^3}}\omega (z+Q_{b\delta})\le (8\theta)^{-1}
\hskip.2in a.s.
\end{align}
and their opposites.

Comparing to the extreme value problem in
the non-critical setting, the strong laws in (\ref{intro-14})
and (\ref{intro-15}) are much more 
sensitive to truncation radius $R$ (more precisely, to the number 
of the $\delta$-neighborhoods that covers the ball $\{\vert x\vert\le R\}$),
as they corresponds to the polynomial (rather than exponential)
decay of the Poissonian tail. To validate the first step
in (\ref{intro-13}) in the argument for the upper bounds, on the other
hand, one has to take
$R$ significantly larger than it is in the proof for the lower bounds.
The impact of larger $R$ in (\ref{intro-14})
and (\ref{intro-15}) can be counter-balanced by taking smaller $\delta$,
even though this action leads to a further increase of  the number
of the $\delta$-neighborhoods.
The cost of this strategy turns out to be a possibly very large
function $A(t)$ appearing on the right hand side of (\ref{intro-13}). An 
observation unique to the critical setting is the irrelevance of
$A(t)$ to
the asymptotic behaviors of the system, as the quantity
$$
H\Big(\theta\max_{\stackrel{\scriptstyle\vert z\vert\le R}
{\scriptstyle z\in 2\delta\Z^3}}\omega (z+Q_{b\delta})\Big)
$$
is equal to zero eventually (under (\ref{intro-14})) or infinitely often
(under (\ref{intro-15})).

\bigskip

\section{Basic estimates}\label{estimate}

In this section we give some auxiliary results that 
will be used in our proofs. We state them
separately for a convenient reference. For future reference, all results in
this section are established in the space $\R^d$ for 
$d\ge 1$, except Lemma \ref{FK-1}
where $d=3$.

\subsection{Association of infinitely divisible fields}\label{sec:asso}

Recall that random variables $X_1,\dots,X_n$ are said to be {\em associated} if for any bounded measurable functions $f, g: \R^n \mapsto \R$ 
non-decreasing (equivalently, non-increasing) in each coordinate
\begin{equation} \label{eq:asso}
\mathrm{Cov}(f(X), g(X)) \ge 0.
\end{equation}
Association is a fairly strong property exhibiting positive dependence. 

Consider now a non-negative random measure $M$ on $\R^d$, taking independent values on disjoint sets such that $M(A)$ is infinitely divisible with the characteristic function
\begin{equation} \label{eq:id-ch}
\E \exp\{- u M(A)\} = \exp\Big\{-m(A) \int_0^{\infty} (1- e^{-us}) \, \rho(ds)\Big\}, \quad u>0,
\end{equation}
for every Borel set $A \subset \R^d$ with $m(A) < \infty$,
where $m$ is a $\sigma$-finite measure on $\R^d$ and $\rho$ is a measure on $(0,\infty)$ such that ${\int_0^{\infty} \min\{s,1\} \, \rho(ds) <\infty}$.  $M$ can be viewed as a  distribution of obstacles in $\R^d$ having random locations and random masses, so we call it an {\em infinitely divisible random field}. $M$ is a Poisson field if $m(dx)=dx$ and $\rho(ds)=\delta_1(ds)$. See \cite{R-R} for more information on infinitely divisible random measures. 

\begin{lemma}\label{lem:00}
Let $M$ be infinitely divisible random field. Then  $M(A_1),\dots,M(A_n)$ are associated for any Borel sets $A_j$ with $m(A_j) < \infty$, $j=1,\dots,n$.  In particular, for all $c_1,\dots,c_n \in \R$
\begin{equation} \label{eq:a1}
\P( M(A_1) \le c_1, \dots, M(A_n) \le c_n) \ge \prod_{j=1}^n 
\P( M(A_j) \le c_j)
\end{equation}
and
\begin{equation} \label{eq:a2}
\P( M(A_1) \ge c_1, \dots, M(A_n) \ge c_n) \ge \prod_{j=1}^n \P( M(A_j) \ge c_j). 
\end{equation}

\end{lemma}

\proof
It follows from \eqref{eq:id-ch} that 
$X=(M(A_1),\dots,M(A_n))$ has infinitely 
divisible distribution without  Gaussian part and its L\'evy measure is 
concentrated on $\R^d_{+}$. Thus the components of $X$ are 
associated \cite{Res} (see also \cite{Sam} for more 
information on association of infinitely divisible random vectors).

Applying \eqref{eq:asso} recursively for $f=\prod_{j=1}^{n-1} \1_{(-\infty, c_j]}$ and $g=\1_{(-\infty, c_n]}$ ($f=\prod_{j=1}^{n-1} \1_{[c_j, \infty)} $ and $g=\1_{[c_n,\infty)}$, respectively) we obtain \eqref{eq:a1} (\eqref{eq:a2}, respectively). 
\qed

\subsection{Truncating Poisson potentials}\label{poisson-1}

In this subsection we study a family of Poisson potentials generated by a smooth truncation of the singular potential kernel $K(x)=|x|^{-p}$, where  $p \in (d/2,d)$. 
The following notation will be used throughout this paper:

$\alpha$: $\R^+\longrightarrow [0,1]$ denotes a fixed smooth function
with the following properties: $\alpha(\lambda)=1$ on $[0,1]$,
$\alpha(\lambda)=0$ for $\lambda\ge 3$ and $-1\le\alpha'(\lambda)\le 0$.
Let $a>0$ be fixed but arbitrary. Put
\begin{align}\label{poisson-2'}
L_{a}(x)={1-\alpha(a^{-1} |x|) \over |x|^p}
\hskip.1in\hbox{and}\hskip.1in
\ol{V}_{a,\epsilon}(x)=\int_{\R^d}L_{a}(x-y)\big[\omega(\epsilon dy)
-\epsilon dy\big].
\end{align}

Let $D\subset\R^d$ be a fixed bounded set.

\begin{lemma}\label{poisson-2} For any $\theta>0$ and fixed $a>0$
\begin{align}\label{poisson-3}
\E\exp\Big\{\theta\sup_{x\in D}
\vert\ol{V}_{a,1}(x)\vert\Big\}
<\infty.
\end{align}

Further, for given $\theta>0$ one can take $a>0$ large enough so
\begin{align}\label{poisson-4}
\sup_{0<\epsilon <1}
\E\exp\Big\{\theta(\log\epsilon^{-1})\sup_{x\in D}
\vert\ol{V}_{a,\epsilon}(x)\vert\Big\}
<\infty.
\end{align}
\end{lemma}

\proof Due to similarity, we only prove (\ref{poisson-4}).
Write
$$
\Psi(\lambda)=e^\lambda-1-\lambda\hskip.2in \lambda\in\R.
$$
We have for $\theta>0$
$$
\begin{aligned}
\E\exp\Big\{\pm\theta(\log\epsilon^{-1})\ol{V}_{a,\epsilon}(0)\Big\}
&=\exp\bigg\{\epsilon\int_{\R^d}\Psi\Big(\pm\theta(\log\epsilon^{-1}) L_a(y)\Big) \, dy\bigg\}\\
&\le \exp\bigg\{\epsilon\int_{\R^d}\Psi\Big(\theta(\log\epsilon^{-1}) L_a(y)\Big) \, dy\bigg\},
\end{aligned}
$$
where the inequality follows from the fact that 
$\Psi(-\lambda)\le\Psi(\lambda)$ for any $\lambda \ge 0$.
Therefore,
$$
\E\exp\Big\{\theta(\log\epsilon^{-1})\vert\ol{V}_{a,\epsilon}(0)\vert\Big\}
\le 2\exp\bigg\{\epsilon\int_{\R^d}\Psi\Big(\theta(\log\epsilon^{-1}) L_a(y)\Big) \, dy\bigg\}.
$$

By a change of variable,
\begin{align} \label{poisson-5}
\int_{\R^d} & \Psi\Big(\theta(\log\epsilon^{-1}) L_a(y)\Big) \, dy 
=(\log\epsilon^{-1})^{d/p}\int_{\R^d}\Psi\Big(\theta L_{a(\log\epsilon^{-1})^{-1/p}}(x) \Big) \, dx\nonumber\\
&\le (\log\epsilon^{-1})^{d/p}
\int_{\{\vert x\vert\ge a(\log\epsilon^{-1})^{-1/p}\}}
\Psi\Big(\theta |x|^{-p}\Big) \, dx\nonumber\\
& \le (\log\epsilon^{-1})^{d/p}\bigg\{
\int_{\{1 \wedge a (\log\epsilon^{-1})^{-1/p}\le \vert x\vert\le 1\}}
+\int_{\{\vert x\vert\ge 1\}}\bigg\}\Psi\Big(\theta |x|^{-p} \Big) \, dx \nonumber\\
& \le (\log\epsilon^{-1})^{d/p}\bigg\{C\exp\Big\{\theta a^{-p}\log \epsilon^{-1}\Big\}
+\int_{\{\vert x\vert\ge 1\}}\Psi\Big(\theta |x|^{-p}\Big) \, dx \bigg\}, 
\end{align}
where the last step follows  from the bound 
$\Psi(\lambda)\le e^{\lambda}$ ($\lambda>0$).
Since
$$
\int_{\{\vert x\vert\ge 1\}}\Psi\Big(\theta |x|^{-p}\Big)dx<\infty,
$$
we get for $a > \theta^{1/p}$ 
\begin{align}\label{poisson-6}
\sup_{0<\epsilon<1} 
\E\exp\Big\{\theta(\log\epsilon^{-1})\vert\ol{V}_{a,\epsilon}(0)\vert\Big\} <\infty.
\end{align}

Similarly as at the beginning of the proof, for any $x, y\in D$ with $x\not =y$,
$$
\begin{aligned}
\E & \exp\bigg\{\theta(\log\epsilon^{-1}){\vert\ol{V}_{a,\epsilon}(x)-
\ol{V}_{a,\epsilon}(y)\vert\over\vert x-y\vert}\bigg\}\\
&\le 2\exp\bigg\{\epsilon\int_{\R^d}\Psi\bigg({\theta\log\epsilon^{-1}\over
\vert x-y\vert}\vert L_a(x-z)-L_a(y-z)\vert\bigg) \, dz\bigg\}.
\end{aligned}
$$
By the mean value theorem we obtain for all $x, y \in \R^d$
\begin{align*}
|L_a(z-x)-L_a(z-y)| & \le C a^{-1} |x-z|^{-p}\1(|x-z|>a) |x-y|\\
& \quad + C a^{-1} |y-z|^{-p}\1(|y-z|>a) |x-y|,
\end{align*}
where $C=p+1$.
Using this estimate and the convexity of $\Psi$ we get
\begin{align*}
\int_{\R^d} & \Psi\bigg({\theta\log\epsilon^{-1}\over
\vert x-y\vert}\vert L_a(x-z)-L_a(y-z)\vert\bigg) \, dz\\
&\le \frac{1}{2} \int_{\{|x-z|>a\}}
\Psi\Big(2\theta Ca^{-1}\log\epsilon^{-1} |x-z|^{-p}) \Big) \, dz\\
& \quad + \frac{1}{2} \int_{\{|y-z|>a\}}
\Psi\Big(2\theta Ca^{-1}\log\epsilon^{-1} |y-z|^{-p}) \Big) \, dz\\
&=(\log \epsilon^{-1})^{d/p}\int_{\{\vert z\vert\ge a
(\log\epsilon^{-1})^{-1/p}\}}
\Psi\Big(2\theta Ca^{-1}|z|^{-p}\Big) \, dz.
\end{align*}

By the same estimate as in (\ref{poisson-5}) we get
$$
\sup_{0<\epsilon<1}\sup_{x\not =y}\E\exp
\bigg\{\theta(\log\epsilon^{-1}){\vert\ol{V}_{a,\epsilon}(x)-
\ol{V}_{a,\epsilon}(y)\vert\over\vert x-y\vert}\bigg\}<\infty.
$$
By Theorem D.6, p.313, \cite{Chen} we 
\begin{align}\label{poisson-7}
\sup_{0<\epsilon<1}\E
\exp\bigg\{\theta(\log\epsilon^{-1})\sup_{x,y\in D}\vert\ol{V}_{a,\epsilon}(x)-
\ol{V}_{a,\epsilon}(y)\vert\bigg\}<\infty.
\end{align}
So the desired conclusion follows from (\ref{poisson-6}) 
and (\ref{poisson-7}).
\qed

Using above lemma, we derive the following almost sure bounds

\begin{lemma}\label{poisson-8} For any $a>0$
\begin{align}\label{poisson-9}
\lim_{R\to\infty}(\log R)^{-1}\sup_{\vert x\vert\le R}\vert\ol{V}_{a, 1}(x)\vert
=0\hskip.2in a.s.
\end{align}

Further, for any positive sequence $\epsilon_n$ such that
$$
\limsup_{n\to\infty}{\epsilon_{n+1}\over\epsilon_n}<1
$$
and any constants $\beta>0$, the strong law
\begin{align}\label{poisson-10}
\lim_{n\to\infty}\sup_{\vert x\vert\le\epsilon_n^{-\beta}}
\vert\ol{V}_{a, \epsilon_n}(x)\vert
=0\hskip.2in a.s.
\end{align}
holds for sufficiently large $n$.
\end{lemma}

\proof The ball $\{x\in\R^d;\hskip.05in\vert x\vert\le R\}$ 
is covered by roughly $CR^d$ unit balls.
By homogeneity of the field $\ol{V}(\cdot)$, for each $\delta>0$
$$
\begin{aligned}
&\P\Big\{\sup_{\vert x\vert\le R}\vert\ol{V}_{a, 1}(x)\vert\ge\delta\log R\Big\}
\le CR^d\P\Big\{\sup_{\vert x\vert\le 1}\vert\ol{V}_{a, 1}(x)\vert
\ge\delta\log R\Big\}\\
&\le CR^{-(\theta\delta -d)}\E\exp\Big\{\theta\sup_{\vert x\vert\le 1}
\vert\ol{V}_{a, 1}(x)\vert\Big\}.
\end{aligned}
$$
Take $\theta$ sufficiently large so $\theta\delta -d\ge 1$.
By (\ref{poisson-3}) the exponential moment on the right hand side
is finite. Thus,
$$
\sum_n\P\Big\{\sup_{\vert x\vert\le 2^n}\vert\ol{V}_{a, 1}(x)
\vert\ge\delta\log 2^n\Big\}<\infty.
$$
Notice that $\delta>0$ can be arbitrarily small.
By Borel-Cantelli lemma
$$
\lim_{n\to\infty}(\log 2^n)^{-1}\sup_{\vert x\vert\le 2^n}\vert\ol{V}_{a, 1}(x)
\vert =0\hskip.2in a.s.
$$
Hence, (\ref{poisson-9}) follows from the fact that
$\displaystyle\sup_{\vert x\vert\le R}\vert\ol{V}_{a, 1}(x)\vert$
is non-decreasing in $R$.

We now come to the proof of (\ref{poisson-10}).
First notice that the ball $B(0, \epsilon_n^{-\beta})$ can be covered by 
$C\epsilon_n^{-d\beta}$ balls of radius 1.
Thus, for each $\delta>0$
$$
\begin{aligned}
&\P\Big\{\sup_{\vert x\vert\le\epsilon_n^{-\beta}}\vert 
\ol{V}_{a,\epsilon_n}(x)\vert\ge \delta\Big\}
\le C\epsilon_n^{-d\beta}\P\{\sup_{\vert x\vert\le 1}\vert 
\ol{V}_{a,\epsilon_n}(x)\vert\ge \delta\Big\}\\
&\le C \epsilon_n^{\theta\delta -d\beta}
\E\exp\Big\{\theta\big(\log \epsilon_n^{-1}\big)\sup_{\vert x\vert\le 1}
\vert \ol{V}_{a,\epsilon_n}(x)\vert\Big\}.
\end{aligned}
$$
Take $\theta>0$ sufficiently large so $\theta\delta -d\beta\ge 1$. 
By (\ref{poisson-4}),
the exponential moment on the right hand side is bounded
uniformly over $n$ when $a>0$ is sufficiently large. Hence,
$$
\sum_n\P\Big\{\sup_{\vert x\vert\le \epsilon_n^{-\beta}}\vert 
\ol{V}_{a,\epsilon_n}(x)\vert\ge \delta\Big\}<\infty.
$$
Therefore, (\ref{poisson-10}) follows from Borel-Cantelli lemma. \qed

\subsection{Lower bound on Brownian motion before it exits a ball}

$B(x,r)$ will denote a  ball in $\R^d$ with center at $x$ and radius $r$.

\begin{lemma}\label{exit}
For every $R, t>0$ and a Borel set $A \subset \R^d$
\begin{align} \label{eq:lb-box}
\P_0\big(B_t \in A, & \,  \max_{0\le s\le t} |B_s| \le 2R\big) \\
& \ge \P_0\big(B_t \in A\cap B(0, R) \big) \ \P_0\big(\max_{0\le s\le 1} |B^0_s| \le Rt^{-1/2}\big),\nonumber
\end{align}
where $B^0_s$ is the Brownian bridge in $\R^d$.  
\end{lemma}

\proof 
Let $B^0_s= B_s-sB_1$ be a Brownian bridge on $[0,1]$ taking values in $\R^d$ and let $Z$ be independent of $B_s$ standard Gaussian random vector in $\R^d$. Then $X_s=B_s^0 + sZ$ is a standard Brownian motion on $0\le s \le 1$ in $\R^d$. 
We have by scaling
\begin{align*}
\P(B_t \in A, \, & \max_{0\le s\le t} |B_s| \le 2R) = \P(t^{1/2} B_1 \in A, \, \max_{0\le s\le 1} |B_s| \le 2Rt^{-1/2}) \\ & = \P(t^{1/2} Z \in A, \, \max_{0\le s\le 1} |B^0_s+ sZ| \le 2Rt^{-1/2}) \\
& \ge \P(t^{1/2} Z \in A, \, |Z| \le Rt^{-1/2}, \, \max_{0\le s\le 1} |B^0_s| \le Rt^{-1/2}) \\
& = \P(B_t \in A\cap B(0,R))\, \P(\max_{0\le s\le 1} |B^0_s| \le Rt^{-1/2}), 
\end{align*}
which proves \eqref{eq:lb-box}.
\qed

\subsection{Bounds by Feynman-Kac functionals}\label{FK}

Given a bounded open domain $D\subset\R^d$, let $W^{1,2}(D)$ be 
 the Sobolev space over $D$, defined to be the 
closure of the inner product space
consists of the infinitely differentiable functions compactly supported in
$D$ under the Sobolev norm
$$
\|g\|_H=\Big\{\|g\|^2_{{\cal L}^2(D)}+\|\nabla g\|^2_{{\cal L}^2(D)}\Big\}^{1/2}.
$$
Define
\begin{align}\label{FK-0}
{\cal F}_d(D)=\bigg\{g\in W^{1,2}(D);\hskip.1in \int_D g^2(x)dx=1 \bigg\}.
\end{align}
For any measurable function $\zeta$ on $D$, put
\begin{align}\label{FK-0'}
\lambda_\zeta(D)=\sup_{g\in{\cal F}_d(D)}\bigg\{\int_D\zeta(x)g^2(x)dx-{1\over 2}
\int_D\vert\nabla g(x)\vert^2dx\bigg\}.
\end{align}
Clearly, $\lambda_\xi(D)\le \lambda_\eta(D)$ and   
$\lambda_\zeta(D)\le \lambda_\zeta(D')$ whenever $\xi(x)\le \eta(x)$
($x\in D$) and $D\subset D'$.

Let
$$
\tau_D=\inf\{s\ge 0;\hskip.1in B_s\not\in D\}.
$$
For any $r>0$, define $T_r=\tau_{B(0,r)}$ with 
$B(0,r)=\{x\in\R^d;\hskip.1in\vert x\vert <r\}$.

\begin{lemma}\label{FK-6} Let $R>0$ and let
$\zeta(x)$ be a function on $\R^d$ such that
$K=\sup_{x\in B(0, 2R)}\zeta(x)<\infty$. We have that for any $t, t_0>0$
satisfying $t_0<t$,
\begin{align}\label{FK-6-1}
\int_{B(0, R)} & \E_x\Bigg[\exp\bigg\{\int_0^t\zeta(B_s)ds\bigg\}; 
T_{R}\ge t\Bigg]dx\\
&\ge (2\pi t_0)^{d/2}e^{-t_0K}\exp\Big\{(t+t_0)
\lambda_\zeta \big(B(0,R)\big)\Big\},\nonumber
\end{align}
\begin{align}\label{FK-6-2}
\E_0 & \Bigg[\exp\bigg\{\int_0^t\zeta(B_s)ds\bigg\}; T_{2R}\ge t\Bigg]\\
&\ge \P_0\Big\{\max_{0\le s\le 1}\vert B_s^0\vert\le Rt_0^{-1/2}\Big\}
\exp\Big\{-2t_0K-{R^2\over 2t_0}\Big\}
\exp\Big\{-t\lambda_{\zeta}\big(B(0,R)\big)\Big\}\nonumber
\end{align}
where $B_s^0$ ($0\le s\le 1$) is a $d$-dimensional Brownian bridge.
\end{lemma}

\proof By a standard procedure of
approximation we may assume that $\zeta(\cdot)$ is H\"older continuous.
 By Feynman-Kac representation
$$
u(t,x)=\E_x\Bigg[
\exp\bigg\{\int_0^{t} \zeta(B_s)ds\bigg\}dx;\hskip.05in
T_R\ge t\Bigg]
$$
solves the initial-boundary value problem
$$
\begin{aligned}
\left\{\begin{array}{ll}\displaystyle
\partial_tu(t,x)={1\over 2}\Delta u(t,x)+\zeta(x)u(t,x)
\hskip.2in (t,x)\in (0, \infty)
\times B(0,R),\\\\
u(0, x)=1\hskip 1.6in  x\in B(0, R) ,\\\\
u(t, x)=0\hskip 1.6in t> 0 \hskip.05in\hbox{and }\hskip.05in\vert x\vert=R.\end{array}\right.
\end{aligned}
$$

Let $\lambda_1>\lambda_2\ge\lambda_3\ge\cdots$ be
the eigenvalues of the operator $(1/2)\Delta+\zeta$ in 
${\cal L}^2\big(B(0, R)\big)$
with zero boundary condition and initial value 1 in $B(0, R)$ and let 
$e_k\in {\cal L}^2\big(B(0, R)\big)$
be an orthogonal basis corresponding to $\{\lambda_k\}$. By (2.31)
in \cite{GKM},
$$
\E_x\exp\bigg\{\int_0^{t} \zeta(B_s)ds\bigg\}\delta_x(B_t);\hskip.05in
T_R\ge t\Bigg]=\sum_{k=1}^\infty e^{t\lambda_k} e^2_k(x)
\ge e^{t\lambda_1} e^2_1(x)
$$
where $\delta_x(\cdot)$ is the Dirac function on $\R^d$
with concentration at $x$.

Noticing the fact that $\lambda_1=\lambda_\zeta\big(B(0, R)\big)$ 
and integrating both sides
we have
$$
\int_{B(0, R)}\E_x\exp\bigg\{\int_0^{t} \zeta(B_s)ds\bigg\}
\delta_x(B_t);\hskip.05in
T_R\ge t\Bigg]dx\ge \exp\Big\{t\lambda_\zeta\big(B(0, R)\big)\Big\}.
$$

In addition, by Markov property
$$
\begin{aligned}
\E_x & \Bigg[\exp\bigg\{\int_0^{t}\zeta(B_s)ds\bigg\}
\delta_x(B_{t+1});\hskip.05in
T_R\ge t\Bigg]\\
&\le e^{t_0K}\E_x\Bigg[\exp\bigg\{\int_0^{t-t_0}\zeta(B_s)ds\bigg\}
\delta_x(B_{t});\hskip.05in
T_R\ge t-t_0\Bigg]\\
&=e^{t_0K}\E_x\Bigg[\exp\bigg\{\int_0^{t-t_0}\zeta(B_s)ds\bigg\}
p_{t_0}(B_{t-t_0}-x);\hskip.05in
T_R\ge t-t_0\Bigg]
\end{aligned}
$$
where
$$
p_{t_0}(y)={1\over (2\pi t_0)^{d/2}}\exp\Big\{-{\vert y\vert^2\over 2t_0}\Big\}\le
{1\over (2\pi t_0)^{d/2}}\hskip.2in y\in\R^d.
$$
Hence, we have proved that
\begin{align}\label{FK-1'}
\int_{B(0, R)} & \E_x\exp\bigg\{\int_0^{t-t_0} \zeta(B_s)ds\bigg\};\hskip.05in
T_R\ge t-t_0\Bigg]dx\\
&\ge (2\pi t_0)^{d/2}e^{-t_0K}
\exp\Big\{t\lambda_\zeta\big(B(0, R)\big)\Big\}.\nonumber
\end{align}
Replacing $t$ by $t+t_0$ leads to (\ref{FK-6-1}).

On the other hand, using Markov property again
$$
\begin{aligned}
\E_0 & \Bigg[\exp\bigg\{\int_0^t\zeta(B_s)ds\bigg\}; T_{2R}\ge t\Bigg]\ge
e^{-t_0K}\E_0\Bigg[\exp\bigg\{\int_{t_0}^t\zeta(B_s)ds\bigg\}; T_{2R}\ge t\Bigg]\\
&=e^{-t_0K}\int_{B(0,2R)}\tilde{p}_{t_0}(x)\E_x\Bigg[\exp\bigg\{\int_0^{t-t_0}
\zeta(B_s)ds\bigg\}; T_{2R}\ge t-t_0\Bigg]dx\\
&\ge e^{-t_0K}\int_{B(0,R)}\tilde{p}_{t_0}(x)\E_x\Bigg[\exp\bigg\{\int_0^{t-t_0}
\zeta(B_s)ds\bigg\}; T_{R}\ge t-t_0\Bigg]dx 
\end{aligned}
$$
where $\tilde{p}_{t_0}(x)$ is the density function of the measure
$$
\mu_{t_0}(A)=\P_0\Big\{B_{t_0}\in A,\hskip.05in T_{2R}\ge t_0\Big\}
\hskip.2in A\subset\R^d.
$$
By (\ref{eq:lb-box}),
$$
\begin{aligned}
\tilde{p}_{t_0}(x)&\ge\P_0\Big\{\max_{0\le s\le 1}\vert B_s^0\vert\le 
Rt_0^{-1/2}\Big\}
{1\over (2\pi t_0)^{d/2}}\exp\Big\{-{\vert x\vert^2\over 2t_0}\Big\}\\
&\ge \P_0\Big\{\max_{0\le s\le 1}\vert B_s^0\vert\le Rt_0^{-1/2}\Big\}
{1\over (2\pi t_0)^{d/2}}\exp\Big\{-{R^2\over 2 t_0}\Big\}
\hskip.2in x\in B(0, R).
\end{aligned}
$$
Thus, we conclude that
$$
\begin{aligned}
\E_0 & \Bigg[\exp\bigg\{\int_0^t\zeta(B_s)ds\bigg\}; T_{2R}\ge t\Bigg]\\
&\ge
e^{-t_0K}\P_0\Big\{\max_{0\le s\le 1}\vert B_s^0\vert\le Rt_0^{-1/2}\Big\}
{1\over (2\pi t_0)^{d/2}}\exp\Big\{-{R^2\over 2t_0}\Big\}\\
& \quad \times\int_{B(0,R)}\E_x\Bigg[\exp\bigg\{\int_0^{t-t_0}
\zeta(B_s)ds\bigg\}; T_{R}\ge t-t_0\Bigg]dx.
\end{aligned}
$$
Finally, (\ref{FK-6-2})  follows from (\ref{FK-1'}).
\qed

\begin{lemma}\label{FK-1} Let $d=3$. For any $\delta>0$ with 
$\{\vert x\vert\le\delta\}\subset D$,
$$
\begin{aligned}
\E_0 & \Bigg[\exp\bigg\{\theta\int_0^t\ol{V}(B_s)ds\bigg\};\hskip.05in
\tau_D\ge 2t\Bigg]\\
&\le\exp\Big\{\theta t\sup_{\vert x\vert\le\delta/2}\vert\ol{V}_{{\delta\over 6}, 1}
(x)\vert\Big\}+
{6\vert D\vert\over\pi\delta^3}
\E_0\exp\Big\{\sqrt{2\delta}T_1\theta\sup_{x\in D}
\vert\ol{V}_{{\delta\over 6}, 1}
(x)\vert\Big\}
\exp\Big\{t\lambda_{\theta\ol{V}}(D)\Big\}
\end{aligned}
$$
conditioning on the event 
$\Big\{\omega\big\{\vert x\vert\le\delta\big\}=0\Big\}$, where $\vert D\vert$
is the volume of $D$ and $\ol{V}_{a,\epsilon}(\cdot)$ is defined in 
(\ref{poisson-2'}).
\end{lemma}

\proof Notice that $\alpha (\lambda)=0$ for $\lambda\ge 3$. Thus,
on the event $\{\omega(\vert x\vert\le\delta)=0\}$,
$$
\int_{\R^3}{\alpha(6\delta^{-1}\vert x\vert)\over\vert x\vert^2}\omega(dy)
\le\int_{\{\vert y\vert>\delta\}}
{1_{\{\vert y-x\vert\le\delta/2\}}\over\vert x\vert^2}\omega(dy)=0
$$
whenever $\vert x\vert\le\delta/2$. Consequently,
$$
\ol{V}(x)=-C_\delta +\ol{V}_{{\delta\over 6}, 1}(x)\le 
\ol{V}_{{\delta\over 6}, 1}(x)\hskip.2in 
$$
where
$$
C_\delta=\int_{\R^3}
{\alpha(6\delta^{-1}\vert x\vert)\over\vert x\vert^2}dx.
$$

For any $r<\delta/2$, therefore
$$
\int_0^{T_r\wedge t}\ol{V}(B_s)ds=\int_0^{\tau_r\wedge t}
\ol{V}_{{\delta\over 6}, 1}(B_s)ds\le t\sup_{\vert x\vert\le\delta}
\vert\ol{V}_{{\delta\over 6}, 1}(x)\vert.
$$
Thus,
\begin{align}\label{FK-2}
\E_0 & \Bigg[\exp\bigg\{\theta\int_0^t\ol{V}(B_s)ds\bigg\};
\hskip.05in\tau_D\ge 2t
\Bigg]\\
&\le\E_0\Bigg[\exp\bigg\{\theta\int_0^t\ol{V}(B_s)ds\bigg\};
\hskip.05in
T_r\le t, \hskip.05in\tau_D\ge 2t\Bigg]
+\exp\Big\{t\sup_{\vert x\vert\le\delta}\vert\ol{V}_{{\delta\over 6}, 1}(x)\vert\Big\}.
\nonumber
\end{align}

Write $\tau_D'=\inf\{t\ge T_r;\hskip.1in B_s\not\in D\}$.
By Markov property,
$$
\begin{aligned}
\E_0 & \Bigg[\exp\bigg\{\theta\int_0^t\ol{V}(B_s)ds\bigg\};\hskip.05in
T_r\le t,\hskip.05in\tau_D\ge 2t\Bigg]\\
&\le \E_0\Bigg[
\exp\Big\{\theta T_r \Big(\sup_{\vert x\vert\le\delta}
\vert\ol{V}_{{\delta\over 6}, 1}(x)\vert -C_\delta\Big)\Big\}
\exp\bigg\{\theta\int_{T_r}^t\ol{V}(B_s)ds\bigg\};\hskip.05in
T_r\le t,\hskip.05in\tau_D'\ge 2t\Bigg]\\
&=\E_0\Bigg[
\exp\Big\{\theta T_r \Big(\sup_{\vert x\vert\le\delta}\vert\ol{V}_{{\delta\over 6}, 1}
(x)\vert -C_\delta\Big)\Big\}
u_0(t-T_r, B_{T_r});\hskip.05in T_r\le t\Bigg],
\end{aligned}
$$
where 
$$
\begin{aligned}
u_0(s, x) & =\E_x\Bigg[\exp\bigg\{\theta\int_0^s\ol{V}(B_u)du\bigg\};\hskip.05in
\tau_D\ge t+s\Bigg]\\
&\le\E_x\Bigg[\exp\bigg\{\theta\int_0^s\ol{V}(B_u)du\bigg\};\hskip.05in
\tau_D\ge t\Bigg]\equiv u_1(s, x)\hskip.2in (0\le s\le t).
\end{aligned}
$$

Notice that on $\{T_r\le t\}$
$$
\begin{aligned}
u_1(t-T_r, B_{T_r})&\le\exp\Big\{-\theta T_r\inf_{x\in D}\ol{V}(x)\Big\}
\E_{B_{T_r}}\Bigg[\exp\bigg\{\theta\int_0^t\ol{V}(B_u)du\bigg\};\hskip.05in
\tau_D\ge t\Bigg]\\
&=\exp\Big\{-\theta T_r\inf_{x\in D}\ol{V}(x)\Big\}u_2(t, B_{T_r})
\hskip.2in (\hbox{say}).
\end{aligned}
$$

Summarizing our estimate,
$$
\begin{aligned}
\E_0 & \Bigg[\exp\bigg\{\theta\int_0^t\ol{V}(B_s)ds\bigg\};\hskip.05in
T_r\le t,\hskip.05in\tau_D\ge 2t\Bigg]\\
&\le \E_0\Bigg[\exp\bigg\{\theta T_r \Big(\sup_{\vert x\vert\le\delta}
\vert\ol{V}_{{\delta\over 6}, 1}
(x)\vert -C_\delta-\inf_{x\in D}\ol{V}(x)\Big)\bigg\}u_2(t, B_{T_r})\Bigg].
\end{aligned}
$$

Recall the classic facts that $T_r$ and $B_{T_r}$ are independent and that
$B_{T_r}$ is uniformly distributed on the sphere $\{\vert x\vert =r\}$.
So the right hand side is equal to
$$
\E_0\exp\bigg\{\theta T_r \Big(\sup_{\vert x\vert\le\delta}
\vert\ol{V}_{{\delta\over 6}, 1}
(x)\vert -C_\delta-\inf_{x\in D}\ol{V}(x)\Big)\bigg\}
{1\over 4\pi r^2}\int_{\{\vert x\vert =r\}}u_2(t,x)dx.
$$
Using fact that $\{\vert x\vert\le\delta\}\subset D$ and the bound
$$
-\ol{V}(x)\le C_\delta-\ol{V}_{{\delta\over 6}, 1}(x)
\le C_\delta+\vert\ol{V}_{{\delta\over 6}, 1}(x)\vert
$$
we have that
$$
\begin{aligned}
\E_0 & \exp\bigg\{\theta T_r \Big(\sup_{\vert x\vert\le\delta}
\vert\ol{V}_{{\delta\over 6}, 1}
(x)\vert -C_\delta-\inf_{x\in D}\ol{V}(x)\Big)\bigg\}
\le \E_0\exp\Big\{2\theta T_r\sup_{x\in D}
\vert\ol{V}_{{\delta\over 6}, 1}(x)\vert\Big\}\\
&\le \E_0\exp\Big\{2\theta T_\delta\sup_{x\in D}
\vert\ol{V}_{{\delta\over 6}, 1}(x)\vert\Big\}
=\E_0\exp\Big\{\sqrt{2\delta}\theta T_1\sup_{x\in D}
\vert\ol{V}_{{\delta\over 6}, 1}(x)\vert\Big\}.
\end{aligned}
$$
Here we have used the fact that 
$T_r\le T_{\delta/2}\buildrel d\over =\sqrt{\delta/2}T_1$.

By (\ref{FK-2}), we conclude that
$$
\begin{aligned}
(4\pi r^2) & \E_0\Bigg[\exp\bigg\{\theta\int_0^t\ol{V}(B_s)ds\bigg\};
\hskip.05in\tau_D\ge 2t\Bigg]\\
&\le\E_0\exp\Big\{\sqrt{2\delta}\theta T_1\sup_{x\in D}
\vert\ol{V}_{{\delta\over 6}, 1}(x)\vert\Big\}
\int_{\{\vert x\vert =r\}}u_2(t,x)dx \\
& \quad +(4\pi r^2)\exp\Big\{t\sup_{\vert x\vert\le\delta}
\vert\ol{V}_{\delta/6, 1}(x)\vert\Big\}. 
\end{aligned}
$$
Integrating the variable $r$ over $[0, \delta/2]$ on the both sides,
$$
\begin{aligned}
\E_0 & \Bigg[\exp\bigg\{\theta\int_0^t\ol{V}(B_s)ds\bigg\};
\hskip.05in\tau_D\ge 2t\Bigg]\\
&\le {6\over \pi\delta^3}\E_0\exp\Big\{\sqrt{2\delta}\theta T_1\sup_{x\in D}
\vert\ol{V}_{{\delta\over 6}, 1}(x)\vert\Big\}
\int_{\{\vert x\vert \le r\}}u_2(t,x)dx\\
& \quad +\exp\Big\{t\sup_{\vert x\vert\le\delta}
\vert\ol{V}_{{\delta\over 6}, 1}(x)\vert\Big\}.
\end{aligned}
$$
Finally, the desired conclusion follows from the bound
$$
\int_{\{\vert x\vert \le r\}}u_2(t,x)dx\le
\int_Du_2(t,x)dx\le\vert D\vert\exp\Big\{t\lambda_{\theta\ol{V}}(D)\Big\}
$$
where the second step follows from Lemma 4.1 in \cite{Chen-1}. \qed

\section{Lower bounds}\label{lb}

We establish the lower bounds requested by Theorem \ref{intro-4},
Theorem \ref{intro-9} and Theorem \ref{intro-11}.
Let $t$ be either fixed (as in Theorem \ref{intro-4} or increase to
infinity (as in Theorem \ref{intro-9} and Theorem \ref{intro-11}).
Let $\epsilon\to 0$ and $R\to\infty$ either as sequences (when $t$ is fixed)
or as functions of $t$ (when $t\to\infty$). The constraint assumed here is
that $R^2\epsilon^{2/3}t^{-1}\ge c$ eventually for some constant $c>0$.
Other relations among the parameters
introduced above will be specified
later according to the context.

By Brownian scaling,
$$
\begin{aligned}
\E_0 & \exp\bigg\{\theta\int_0^{t}\ol{V}(B_s)ds\bigg\}
=\E\exp\bigg\{\theta\int_0^{t\epsilon^{-2/3}}\ol{V}_\epsilon(B_s)ds\bigg\}\\
&\ge \E_0\Bigg[\exp\bigg\{
\theta\int_0^{t\epsilon^{-2/3}}\ol{V}_\epsilon(B_s)ds\bigg\};\hskip.05in
T_{2R}\ge t\epsilon^{-2/3}\Bigg], \nonumber
\end{aligned}
$$
where
$$
\ol{V}_\epsilon(x)=\int_{\R^3}{1\over \vert y-x\vert^2}
\big[\omega(\epsilon dy)-\epsilon dy\big].
$$

Let $r>0$ and $a>0$ be two large but fixed numbers with
$r<a$. Consider the decomposition
$$
\ol{V}_\epsilon(x)=\ol{V}_{a,\epsilon}(x)+V_{a,\epsilon}(x)
-\epsilon C_{a}
$$
where $\ol{V}_{a,\epsilon}(x)$ is defined as in (\ref{poisson-2'}),
$$
V_{a,\epsilon}(x)=\int_{\R^3}{\alpha(a^{-1}\vert y-x\vert)\over\vert y-x\vert^2} \, 
\omega(\epsilon dy)\hskip.1in
\hbox{and}\hskip.1in
C_{a}=\int_{\R^3}{\alpha(a^{-1}\vert x\vert)\over \vert x\vert^2} \, dx.
$$
We have
\begin{align}\label{lb-1}
\E_0 & \exp\bigg\{\theta\int_0^{t}\ol{V}(B_s)ds\bigg\}\\
&\ge 
\exp\Big\{-\theta t\epsilon^{-2/3}\Big(C_{a}\epsilon +
\sup_{x\in B(0, 2R)}\vert 
\ol{V}_{a,\epsilon}(x)\vert\Big)\Big\}\nonumber\\
& \quad \times\E_0\Bigg[\exp\bigg\{
\theta\int_0^{t\epsilon^{-2/3}}V_{a,\epsilon}(B_s)ds\bigg\};\hskip.05in
T_{2R}\ge t\epsilon^{-2/3}\Bigg].\nonumber
\end{align}

Let $\delta>0$ be a small but fixed number satisfying $r+\delta<a$.
For any $z\in 2r\Z^3\cap B(0, R-r)$ and $x\in\R^d$
$$
\begin{aligned}
\theta V_{a,\epsilon}(x) & \ge \theta
\int_{\{\vert y-z\vert\le\delta\}}{\alpha(a^{-1}\vert y-x\vert)
\over\vert y-x\vert^2}
\omega(\epsilon dy)\\
&\ge\theta\omega\big(B(\epsilon^{1/3}z,\epsilon^{1/3}\delta)\big)
{\alpha\big(a^{-1}(\vert z-x\vert+\delta)\big)
\over (\vert z-x\vert +\delta)^2}\equiv\zeta^z_\epsilon(x).
\end{aligned}
$$
Consequently,
$$
\begin{aligned}
\E_0 & \Bigg[\exp\bigg\{
\theta\int_0^{t\epsilon^{-2/3}}V_{a,\epsilon}(B_s)ds\bigg\};\hskip.05in
T_{2R}\ge t\epsilon^{-2/3}\Bigg]\\
&\ge\E_0\Bigg[\exp\bigg\{
\theta\int_0^{t\epsilon^{-2/3}}\zeta^z_\epsilon(B_s)ds\bigg\};\hskip.05in
T_{2R}\ge t\epsilon^{-2/3}\Bigg]\\
&\ge\P_0\Big\{\max_{0\le s\le 1}\vert 
B_s^0\vert\le {R\epsilon^{1/3}\over\delta\sqrt{t}}\Big\}
\exp\Big\{t\epsilon^{-2/3}\lambda_{\zeta^z_\epsilon}\big(B(0, R)\big)\Big\}\\
& \quad \times\exp\Big\{-2\theta\omega\big(B(\epsilon^{1/3}z,
\epsilon^{1/3}\delta)\big)
t\epsilon^{-2/3}-{R^2\epsilon^{2/3}\over 2\delta^2t}\Big\},
\end{aligned}
$$
where the last step follows from (\ref{FK-6-2}) in Lemma \ref{FK-6} with 
$t$ being replaced by $t\epsilon^{-2/3}$
and $t_0=\delta^2t\epsilon^{-2/3}$, and from
the observation
$$
0\le\zeta^z_\epsilon(x)\le\delta^{-2}\theta\omega
\big(B(\epsilon^{1/3}z,\epsilon^{1/3}\delta)\big)\hskip.2in x\in\R^d.
$$
Notice that $B(z, r)\subset B(0, R)$ and that $r+\delta<a$ leads to
$$
\zeta^z_\epsilon(x)=\theta\omega\big(B(\epsilon^{1/3}z,\epsilon^{1/3}\delta)\big)
{1\over (\vert z-x\vert +\delta)^2}\hskip.2in x\in B(z, r).
$$

By substitution $g(x)\mapsto g(x-z)$, therefore,
$$
\begin{aligned}
\lambda_{\zeta^z_\epsilon} & \big(B(0, R)\big)\ge
\lambda_{\zeta^z_\epsilon}\big(B(z, r)\big)\\
&=\sup_{g\in {\cal F}_3(B(z,r))}\bigg\{\omega\big(B(\epsilon^{1/3}z, 
\epsilon^{1/3}\delta)\big)
\theta\int_{B(z,r)}{g^2(x)\over (\vert x-z\vert+\delta)^2}dx-{1\over 2}
\int_{B(z,r)}\vert \nabla g(x)\vert^2dx\bigg\}\\
&=\sup_{g\in {\cal F}_3(B(0,r))}\bigg\{\omega\big(B(\epsilon^{1/3}z,
\epsilon^{1/3}\delta)\big)
\theta\int_{B(0,r)}{g^2(x)\over (\vert x\vert+\delta)^2}dx-{1\over 2}
\int_{B(0,r)}\vert \nabla g(x)\vert^2dx\bigg\}\\
&=H_{r,\delta}\Big(\omega\big(B(\epsilon^{1/3}z,
\epsilon^{1/3}\delta)\big)\theta\Big)
\end{aligned}
$$
where the function $H_{r,\delta}(\cdot)$ is defined as
$$
H_{r,\delta}(\theta)=\sup_{g\in {\cal F}_3(B(0,r))}
\bigg\{
\theta\int_{B(0,r)}{g^2(x)\over (\vert x\vert+\delta)^2}dx-{1\over 2}
\int_{B(0,r)}\vert \nabla g(x)\vert^2dx\bigg\}.
$$

Summarizing our estimates since (\ref{lb-1}),
$$
\begin{aligned}
\E_0 & \exp\bigg\{\theta\int_0^{t}\ol{V}(B_s)ds\bigg\}\\
&\ge\P_0\Big\{\max_{0\le s\le 1}\vert 
B_s^0\vert\le {R\epsilon^{1/3}\over\delta\sqrt{t}}\Big\}
\exp\Big\{-{R^2\epsilon^{2/3}\over 2\delta^2t}\Big\}
\exp\Big\{t\epsilon^{-2/3}H_{r,\delta}
\Big(\omega\big(B(\epsilon^{1/3}z,\epsilon^{1/3}\delta)\big)\theta\Big)\Big\}\\
& \quad \times\exp\Big\{-\theta t\epsilon^{-2/3}
\Big(C_a\epsilon +\sup_{x\in B(0, R)}\vert 
\ol{V}_{a,\epsilon}(x)\vert +
2\omega\big(B(\epsilon^{1/3}z,\epsilon^{1/3}\Big)\Big\}.
\end{aligned}
$$
By the assumption that $R^2\epsilon^{2/3}/t$ is eventually bounded from below,
there is a constant $\gamma>0$ such that
$$
\P_0\Big\{\max_{0\le s\le 1}\vert 
B_s^0\vert\le {R\epsilon^{1/3}\over\delta\sqrt{t}}\Big\}\ge\gamma
$$
eventually holds. Taking maximum over $z\in 2r\Z^3\cap B(0,R-r)$,
\begin{align}\label{lb-2}
\E_0 & \exp\bigg\{\theta\int_0^{t}\ol{V}(B_s)ds\bigg\}\\
&\ge\gamma\exp\Big\{-{R^2\epsilon^{2/3}\over 2\delta^2t}\Big\}
\exp\Big\{t\epsilon^{-2/3}H_{r,\delta}
\Big(\max_{z\in 2r\Z^3\cap B(0,R-r)}\omega\big(B(\epsilon^{1/3}z,\epsilon^{1/3}\delta)\big)\theta\Big)\Big\}\nonumber\\
& \quad \times\exp\Big\{-\theta t\epsilon^{-2/3}
\Big(C_a\epsilon +\sup_{x\in B(0, R)}\vert 
\ol{V}_{a,\epsilon}(x)\vert +2\max_{z\in 2r\Z^3\cap B(0,R-r)}
\omega\big(B(\epsilon^{1/3}z,\epsilon^{1/3}\Big)\Big\}.\nonumber
\end{align}

A version of (\ref{lb-2}) is also needed and is derived as follows:
By the Brownian scaling,
$$
\int_{B(0,\epsilon^{1/3}R)}\E_x\exp\bigg\{\int_0^t\ol{V}(B_s)ds\bigg\}dx
=\epsilon\int_{B(0, R)}\E_x\exp\bigg\{\int_0^{t\epsilon^{-2/3}}\ol{V}_\epsilon
(B_s)ds\bigg\}dx.
$$
Following the decomposition of $\ol{V}_\epsilon(\cdot)$
the same way as above and then applying (\ref{FK-6-1}) (instead of 
(\ref{FK-6-2})) with $t_0=1$, we have
\begin{align}\label{lb-2'}
\int_{B(0,\epsilon^{1/3}R)} & \E_x\exp\bigg\{\int_0^t\ol{V}(B_s)ds\bigg\}dx\\
&\ge(2\pi)^{d/2}\epsilon \exp\Big\{-\delta^{-2}\max_{z\in 2r\Z^3\cap B(0,R-r)}
\omega\big(B(\epsilon^{1/3}z,\epsilon^{1/3})\big)\Big\}\nonumber\\
& \quad \times\exp\bigg\{-t\epsilon^{-2/3}\theta\Big(\epsilon C_a 
+\sup_{x\in B(0,R)}\vert\ol{V}_{a,\epsilon}(x)\vert\Big)\bigg\}\nonumber\\
& \quad \times\exp\bigg\{(1+t\epsilon^{-2/3})H_{r,\delta}
\Big(\max_{z\in 2r\Z^3\cap B(0,R-r)}
\omega\big(B(\epsilon^{1/3}z,\epsilon^{1/3}\delta)\big)\theta\Big)
\bigg\}.\nonumber
\end{align}

\subsection{Lower bound for Theorem \ref{intro-4}}

We show that when $\theta>1/16$,
\begin{align}\label{lb-3}
\E_0\exp\bigg\{\theta\int_0^{t}\ol{V}(B_s)ds\bigg\}=\infty\hskip.1in
a.s.\hskip.2in \forall
t>0.
\end{align}

Let $t$ be fixed.  Taking $\epsilon=2^{-3n}$ and 
$R=\delta 2^{2n}$ in (\ref{lb-2}) gives
\begin{align}\label{lb-4}
\E_0 & \exp\bigg\{\theta\int_0^{t}\ol{V}(B_s)ds\bigg\}\\
&\ge\gamma\exp\Big\{-{2^{2n}\over 2t}\Big\}\exp\bigg\{t2^{2n}
H_{r,\delta}\Big(\theta\max_{z\in 2r\Z^3\cap B(0,\delta 2^{2n}-r)}
\omega\big(B(2^{-n}z, 2^{-n}\delta)\big)\Big)\bigg\}\nonumber\\
& \quad \times\exp\bigg\{-\theta t2^{2n}\bigg\{C_a2^{-3n}+
\sup_{x\in B(0, \delta 2^{2n})}\vert 
\ol{V}_{a,2^{-3n}}(x)\vert \nonumber\\
&\hskip2.0in +2\max_{z\in 2r\Z^3\cap B(0,\delta 2^{2n}-r)}
\omega\big(B(2^{-n}z, 2^{-n}\delta)\big)\Big)\bigg\}.\nonumber
\end{align}

By (\ref{poisson-10}), 
\begin{align}\label{lb-5}
\lim_{n\to\infty}\sup_{x\in B(0, \delta 2^{2n})}\vert 
\ol{V}_{a,2^{-3n}}(x)\vert =0\hskip.2in a.s.
\end{align}
when $a>0$ is sufficiently large.

We now prove that
\begin{align}\label{lb-6}
\limsup_{n\to\infty}\max_{z\in 2r\Z^3\cap B(0,\delta 2^{2n}-r)}
\omega\big(B(2^{-n}z,2^{-n}\delta)\big)= 2\hskip.2in a.s.
\end{align}

By homogeneity and increment independence of the Poisson field,
The random variables
$$
\omega\big(B(2^{-n}z, 2^{-n}\delta)\big);
\hskip.2in z\in 2r\Z^3\cap B(0,\delta 2^{2n}-r)
$$
are i.i.d's. Hence,
$$
\begin{aligned}
\P & \Big\{\max_{z\in 2r\Z^3\cap B(0,\delta 2^{2n}-r)}
\omega\big(B(2^{-n}z, 2^{-n}\delta)\big)\ge 3\Big\} \\
& \le \#\big\{2r\Z^3\cap B(0,\delta 2^{2n}-r)\big\}
\P\Big\{\omega\big(B(0,2^{-n}\delta)\big)\ge 3\Big\}\\
&\le C 2^{6n} \Big((2^{-n}\delta)^3\Big)^3=O\Big(2^{-3n}\Big\}.
\end{aligned}
$$
Thus,
$$
\sum_n\P\Big\{\max_{z\in 2r\Z^3\cap B(0,\delta 2^{2n}-r)}
\omega\big(B(2^{-n}z, 2^{-n}\delta)\big)\ge 3\Big\}<\infty.
$$
By Borel-Cantelli lemma and by the fact that the random variable
$$
\displaystyle\max_{z\in 2r\Z^3\cap B(0,\delta 2^{2n}-r)}
\omega\big(B(z,2^{-n}\delta)\big)
$$ 
takes
integer-values,
$$
\limsup_{n\to\infty}\max_{z\in 2r\Z^3\cap B(0,\delta 2^{2n}-r)}
\omega\big(B(2^{-n}z,2^{-n}\delta)\big)\le 2\hskip.2in a.s.
$$

On the other hand, write 
$A_n=B(0,\delta 2^{2n}-r)\setminus B(0,\delta 2^{2(n-1)})$.
$$
\begin{aligned}
\P & \Big\{\max_{z\in 2r\Z^3\cap A_n}
\omega\big(B(2^{-n}z,2^{-n}\delta)\big)\le 1\Big\}\\
&=\bigg(1-\P\Big\{
\omega\big(B(0,2^{-n}\delta)\big)\ge 2\Big\}
\bigg)^{\#\{2r\Z^3\cap A_n\}}\\
&\le \Big(1- c\delta^32^{-6n}\Big)^{\#\{2r\Z^3\cap A_n\}}
\le \exp\{-c_0\delta^3\}.
\end{aligned}
$$
So we have that
$$
\sum_n\P\Big\{\max_{z\in 2r\Z^3\cap A_n}
\omega\big(B(2^{-n}z, 2^{-n}\delta)\big)\ge 2\Big\}=\infty.
$$
Notice that the sequence
$$
\max_{z\in 2r\Z^3\cap A_n}
\omega\big(B(2^{-n}z, 2^{-n}\delta)\big)\hskip.2in n=1,2,\cdots
$$
is an independent sequence. By Borel-Cantelli lemma
$$
\begin{aligned}
&\limsup_{n\to\infty}\max_{z\in 2r\Z^3\cap B(0,\delta 2^{2n}-r)}
\omega\big(B(2^{-n}z, 2^{-n}\delta)\big)\\
& \quad \ge\limsup_{n\to\infty}\max_{z\in 2r\Z^3\cap A_n}
\omega\big(B(2^{-n}z, 2^{-n}\delta)\big)\ge 2\hskip.2in a.s.
\end{aligned}
$$

By the fact that
$\theta>16^{-1}$ and by Lemma \ref{H-2},
$$
\lim_{\stackrel{\delta\to 0^+}{\scriptstyle r\to\infty}}H_{r,\delta}(2\theta)
=\sup_{g\in {\cal F}_3}\bigg\{2\theta\int_{\R^3}{g^2(x)\over\vert x\vert^2}dx
-{1\over 2}\int_{\R^3}\vert\nabla g(x)\vert^2dx\bigg\}=\infty.
$$
Therefore, one can take $\delta$ sufficiently small, and $r$ sufficiently
large, so we have $H_{r,\delta}(2\theta)> 2\theta +2^{-1}t^{-2}$. 
Finally, the requested
(\ref{lb-3}) follows from (\ref{lb-4}), (\ref{lb-5}), (\ref{lb-6}). \qed

\subsection{Lower bound for Theorem \ref{intro-9}}

Recall that $0<\theta<1/16$ and $k=[(8\theta)^{-1}]$. We prove
\begin{align}\label{lb-7}
\limsup_{t\to\infty}t^{-{k+1\over k-1}}l(t)^{-{2\over 3(k-1)}}
\log\E_0\exp\bigg\{\theta\int_0^t\ol{V}(B_s)ds\bigg\}=\infty\hskip.2in a.s.
\end{align}
under the assumption
\begin{align}\label{lb-8}
\int_1^\infty {dt\over t \cdot l(t)}=\infty.
\end{align}

Taking $t_n=2^n$, $\epsilon=\epsilon_n=t_n^{-{3\over k-1}}l(t_n)^{-{1\over k-1}}$
and $R=R_n=\delta t_n^{k+1\over k-1}l(t_n)^{2\over 3(k-1)}$ in (\ref{lb-2}) gives
\begin{align}\label{lb-9}
\E_0 & \exp\bigg\{\theta\int_0^{t_n}\ol{V}(B_s)ds\bigg\}\\
&\ge\gamma\exp\Big\{-{1\over 2} t_n^{k+1\over k-1}l(t_n)^{2\over 3(k-1)}\Big\}
\nonumber\\
& \quad \times\exp\bigg\{t_n^{k+1\over k-1}l(t_n)^{2\over 3(k-1)}
H_{r,\delta}\Big(\theta\max_{z\in 2r\Z^3\cap B(0,R_n-r)}
\omega\big(B(\epsilon_n^{1/3}z, \epsilon_n^{1/3}\delta)\big)\Big)\bigg\}
\nonumber\\
& \quad \times\exp\bigg\{-\theta t_n^{k+1\over k-1}l(t_n)^{2\over 3(k-1)}
\Big(C_a\epsilon_n+\sup_{x\in B(0, R_n)}\vert 
\ol{V}_{a,\epsilon_n}(x)\vert\nonumber\\
&\hskip2.1in +2\max_{z\in 2r\Z^3\cap B(0,R_n-r)}
\omega\big(B(\epsilon_n^{1/3}z, \epsilon_n^{1/3}\delta)\big)\Big)\bigg\}.
\nonumber
\end{align}

By (\ref{poisson-10}),
\begin{align}\label{lb-10}
\lim_{n\to\infty}\sup_{x\in B(0, R_n)}\vert 
\ol{V}_{a,\epsilon_n}(x)\vert=0\hskip.2in a.s.
\end{align}
as $a>0$ is sufficiently large.

In addition, 
$$
\begin{aligned}
\P & \Big\{\max_{z\in 2r\Z^3\cap B(0,R_n-r)}
\omega\big(B(\epsilon_n^{1/3}z,\epsilon_n^{1/3}\delta)\big)\ge k+2\Big\}\\
&\le Ct_n^{3(k+1)\over k-1}l(t_n)^{2\over (k-1)} 
\P\Big\{\omega\big(B(0,\epsilon_n^{1/3}\delta)\big)\ge k+2\Big\}\\
&\le Ct_n^{-{3\over k-1}}l(t_n)^{-{k\over k-1}}.
\end{aligned}
$$
Consequently,
$$
\sum_n\P\Big\{\max_{z\in 2r\Z^3\cap B(0,R_n-r)}
\omega\big(B(\epsilon_n^{1/3}z,\epsilon_n^{1/3}\delta)\big)\ge k+2\Big\}<\infty.
$$
By Borel-Cantelli lemma,
\begin{align}\label{lb-11}
\limsup_{n\to\infty}\max_{z\in 2r\Z^3\cap B(0,R_n-r)}
\omega\big(B(\epsilon_n^{1/3}z,\epsilon_n^{1/3}\delta)\big)\le k+1\hskip.2in a.s.
\end{align}

On the other hand, let $A_n=B(0,R_n-r)\setminus B(0,R_{n-1}-r)$.
$$
\begin{aligned}
\P & \Big\{\max_{z\in 2r\Z^3\cap A_n}
\omega\big(B(\epsilon_n^{1/3}z,\epsilon_n^{1/3}\delta)\big)\le k\Big\}\\
&=\bigg(1-\P\Big\{
\omega\big(B(0,\epsilon_n^{1/3}\delta)\big)\ge k+1\Big\}
\bigg)^{\#\{2r\Z^3\cap A_n\}}.
\end{aligned}
$$
Hence, 
$$
\begin{aligned}
\P & \Big\{\max_{z\in 2r\Z^3\cap A_n}
\omega\big(B(\epsilon_n^{1/3}z,\epsilon_n^{1/3}\delta)\big)\ge k+1\Big\}\\
&\sim \#\{2r\Z^3\cap A_n\}\P\Big\{
\omega\big(B(0,\epsilon_n^{1/3}\delta)\big)\ge k+1\Big\}
\ge c_0 l(t_n)^{-1},
\end{aligned}
$$
where $c_0>0$ is a constant independent of $n$. By (\ref{lb-8}), 
$$
\sum_n l(t_n)^{-1}=\infty.
$$
By Borel-Cantelli lemma,
\begin{align}\label{lb-12}
\limsup_{n\to\infty} & \max_{z\in 2r\Z^3\cap B(0,R_n-r)}
\omega\big(B(\epsilon_n^{1/3}z,\epsilon_n^{1/3}\delta)\big)\\
&\ge
\limsup_{n\to\infty}\max_{z\in 2r\Z^3\cap A_n}
\omega\big(B(\epsilon_n^{1/3}z,\epsilon_n^{1/3}\delta)\big)
\ge k+1\hskip.2in a.s.\nonumber
\end{align}

By (\ref{lb-9}), (\ref{lb-10}),
(\ref{lb-11}), (\ref{lb-12}),
$$
\begin{aligned}
\limsup_{n\to\infty} & \ t_n^{-{k+1\over k-1}}l(t_n)^{-{2\over 3(k-1)}} \, 
\log\E_0\exp\bigg\{\theta\int_0^{t_n}\ol{V}(B_s)ds\bigg\}\\
&\ge H_{r,\delta}\big((k+1)\theta\big)-2\theta (k+1)-2^{-1}
\hskip.2in a.s.
\end{aligned}
$$
Notice that $(k+1)\theta>8^{-1}$. By Lemma \ref{H-2}, letting
$r\to\infty$ and $\delta\to 0^+$ on the right hand side
leads to (\ref{lb-7}). \qed

\subsection{Lower bound for Theorem \ref{intro-11}}

We prove that
\begin{align}\label{lb-13}
\liminf_{t\to\infty} \, t^{-{k+1\over k-1}}l(t)^{{2\over 3(k-1)}}
\log\E_0\exp\bigg\{\theta\int_0^t\ol{V}(B_s)ds\bigg\}=\infty\hskip.2in a.s.
\end{align}
under the assumption that
\begin{align}\label{lb-14}
\int_1^\infty {1\over t }\exp\big\{-c\cdot l(t)\big\}dt<\infty
\hskip.2in \forall c>0.
\end{align}

This time we use (\ref{lb-2'}) instead of (\ref{lb-2}).
Taking $t_n=2^n$, $\epsilon=\epsilon_n=t_n^{-{3\over k-1}}
l(t_n)\big)^{{1\over k-1}}$,
$R=R_n=t_n^{k+1\over k-1}l(t_n)^{-{2\over 3(k-1)}}$ in (\ref{lb-2'}) gives
\begin{align}\label{lb-15}
\int_{B(0,\epsilon_n^{1/3}R_n)} & \E_x\exp\bigg\{\int_0^{t_n}\ol{V}(B_s)ds\bigg\}dx\\
& \hskip-.2in \ge(2\pi)^{3/2}\epsilon_n \exp\Big\{-\delta^{-2}\max_{z\in 2r\Z^3\cap B(0,R_n-r)}
\omega\big(B(\epsilon_n^{1/3}z,\epsilon_n^{1/3})\big)\Big\}\nonumber\\
&  \times\exp\bigg\{- t_n^{k+1\over k-1}l(t_n)^{-{2\over 3(k-1)}}
\theta\Big(\epsilon C_a 
+\sup_{x\in B(0,R_n)}\vert\ol{V}_{a,\epsilon_n}(x)\vert\Big)\bigg\}\nonumber\\
&  \times\exp\bigg\{\Big(1+ t_n^{k+1\over k-1}l(t_n)^{-{2\over 3(k-1)}}\Big)H_{r,\delta}
\Big(\max_{z\in 2r\Z^3\cap B(0,R_n-r)}
\omega\big(B(\epsilon_n^{1/3}z,\epsilon_n^{1/3}\delta)\big)\theta\Big)
\bigg\}.\nonumber
\end{align}

We now show that for any $\delta>0$ and $r>0$,
\begin{align}\label{lb-16}
\liminf_{n\to\infty}\max_{z\in 2r\Z^3\cap B(0,R_n-r)}
\omega\big(B(\epsilon_n^{1/3}z,\epsilon_n^{1/3}\delta)\big)\ge k+1\hskip.2in a.s.
\end{align}

Indeed, by independence
$$
\begin{aligned}
\P & \Big\{\max_{z\in 2r\Z^3\cap B(0,R_n-r)}
\omega\big(B(\epsilon_n^{1/3}z,\epsilon_n^{1/3}\delta)\big)\le k\Big\}\\
&=\bigg(1-\P\Big\{
\omega\big(B(0,\epsilon_n^{1/3}\delta)\big)\ge k+1\Big\}
\bigg)^{\#\{2r\Z^3\cap B(0,R_n-r)\}}.
\end{aligned}
$$
By the fact that
$$
\begin{aligned}
\P & \Big\{\omega\big(B(0,\epsilon_n^{1/3}\delta)\big)\ge k+1\Big\}
\sim {1\over (k+1)!}\Big({4\over 3}\pi \delta^3\epsilon_n\Big)^{k+1}\\
&={1\over (k+1)!}\Big({4\over 3}\pi\delta^3\Big)^{k+1}
t_n^{-{3(k+1)\over k-1}}l(t_n)^{k+1\over k-1},
\end{aligned}
$$
there is a constant $c=c(k,\delta, r)>0$ such that
$$
\P\Big\{\max_{z\in 2r\Z^3\cap B(0,R_n-r)}
\omega\big(B(\epsilon_n^{1/3}z,\epsilon_n^{1/3}\delta)\big)\le k\Big\}
\le\exp\Big\{-c\cdot l(t_n)\Big\}
$$
for large $n$. By (\ref{lb-14}),
$$
\sum_n\exp\Big\{-c\cdot l(t_n)\Big\}<\infty.
$$
Hence, (\ref{lb-16}) follows from Borel-Cantelli lemma.

Notice that (\ref{lb-10}) and (\ref{lb-11}) remain true in this setting.
By (\ref{lb-15}) and (\ref{lb-16}), therefore,
\begin{align}\label{lb-17}
\liminf_{n\to\infty}t_n^{-{k+1\over k-1}}l(t_n)^{2\over 3(k-1)}
\int_{B(0,R_n\epsilon_n^{1/3})}\E_x\exp\bigg\{\int_0^{t_n}\ol{V}(B_s)ds\bigg\}dx
=\infty.
\end{align}

For any large $t>0$, let $n$ be such that $t_n\le t\le t_{n+1}$.
Then
$$
\E_0\exp\bigg\{\theta\int_0^{t}\ol{V}(B_s)ds\bigg\}
\ge\E_0\Bigg[\exp\bigg\{\theta\int_0^{t}\ol{V}(B_s)ds\bigg\};\hskip.05in 
T_{R_n\epsilon_n^{1/3}}\ge t-t_{n-1}\Bigg].
$$

Notice that for any $x\in\R^d$
$$
\begin{aligned}
\ol{V}(x)= &\int_{\R^3}{\alpha(a^{-1}\epsilon_{n}^{-1/3}\vert y-x\vert)
\over\vert y-x\vert^2}\omega(dy)-C_a\epsilon_{n}^{1/3}\\
&+
\int_{\R^3}{1-\alpha(a^{-1}\epsilon_{n}^{-1/3}\vert y-x\vert)
\over\vert y-x\vert^2}\big[\omega(dy)-dy\big]\\
&\ge -C_a\epsilon_{n}^{1/3}+\epsilon_{n}^{-2/3}
\ol{V}_{a,\epsilon_{n}}(\epsilon_{n}^{-1/3}x).
\end{aligned}
$$
So we have
$$
\begin{aligned}
\E_0 & \exp\bigg\{\theta\int_0^{t}\ol{V}(B_s)ds\bigg\}\\
&\ge\exp\bigg\{-\theta t_{n+1}\Big(C_a\epsilon_{n}^{1/3}+\epsilon_{n}^{-2/3}
\sup_{x\in B(0, R_{n})} \vert\ol{V}_{a,\epsilon_{n}}(x)\vert\Big)\bigg\}\\
& \quad \times\E_0\Bigg[\exp\bigg\{\theta\int_{t-t_{n-1}}^{t}\ol{V}(B_s)
ds\bigg\};\hskip.05in 
T_{R_n\epsilon_n^{1/3}}\ge t-t_{n-1}\Bigg].
\end{aligned}
$$

By Markov property,
$$
\begin{aligned}
\E_0 & \Bigg[\exp\bigg\{\theta\int_{t-t_{n-1}}^{t}\ol{V}(B_s)ds\bigg\};\hskip.05in 
T_{R_{n}\epsilon_n^{1/3}}\ge t-t_{n-1}\Bigg]\\
&=\int_{B(0, R_n\epsilon_n^{-1/3})}\tilde{p}_{t-t_{n-1}}(x)
\E_x\exp\bigg\{\theta\int_0^{t_{n-1}}\ol{V}(B_s)ds\bigg\}dx\\
&\ge \int_{B(0, \hskip.02in R_{n-1}\epsilon_{n-1}^{-1/3})}\tilde{p}_{t-t_{n-1}}(x)
\E_x\exp\bigg\{\theta\int_0^{t_{n-1}}\ol{V}(B_s)ds\bigg\}dx
\end{aligned}
$$
where $\tilde{p}_{t-t_{n-1}}(x)$ is the density function of the measure
$$
\mu_{t-t_{n-1}}(A)=\P_0\big\{B_{t-t_{n-1}}\in A,\hskip.05in 
T_{R_n\epsilon_n^{1/3}}\ge t-t_{n-1}\big\}\hskip.2in A\subset\R^d.
$$
Notice that $R_{n-1}\epsilon_{n-1}^{-1/3}\le 2^{-1} R_n\epsilon_n^{-1/3}$
for large $n$. By Lemma \ref{exit},
$$
\begin{aligned}
&\tilde{p}_{t-t_{n-1}}(x)\\
&\ge \P_0
\Big\{\max_{0\le s\le 1}\vert B_s^0\vert\le R_n\epsilon_n^{1/3}
(t-t_{n-1})^{-1/2}\Big\}\big(2\pi (t-t_{n-1})\big)^{-3/2}
\exp\Big\{-{\vert x\vert^2\over 2(t-t_{n-1})}\Big\}\\
&\ge \P_0\Big\{\max_{0\le s\le 1}\vert B_s^0\vert\le R_n\epsilon_n^{1/3}
t_{n+1}^{-1/2}\Big\}\big(2\pi t_{n+1}\big)^{-3/2}
\exp\Big\{-{R_{n-1}^2\epsilon_{n-1}^{2/3}\over 2t_{n-1}}\Big\}\\
&\ge\gamma 2^{-3n/2}\exp\Big\{-{1\over 2}t_{n-1}^{k+1\over k-1}
l(t_{n-1})^{2\over 3(k-1)}\Big\}
\end{aligned}
$$
for every $x\in B(0, \hskip.02in R_{n-1}\epsilon_{n-1}^{-1/3})$, where
$\gamma>0$ is a constant independent of $t$.

Summarizing our computation,
$$
\begin{aligned}
&\E_0\exp\bigg\{\theta\int_0^{t}\ol{V}(B_s)ds\bigg\}\\
&\ge\gamma 2^{-3n/2}
\exp\bigg\{-\theta t_{n+1}\Big(C_a\epsilon_{n}^{1/3}+\epsilon_{n}^{-2/3}
\sup_{x\in B(0, R_{n})} \vert\ol{V}_{a,\epsilon_{n}}(x)\vert\Big)\bigg\}\\
& \ \times\exp\Big\{-{1\over 2}t_{n-1}^{k+1\over k-1}
l(t_{n-1})^{2\over 3(k-1)}\Big\}
\int_{B(0, \hskip.02in R_{n-1}\epsilon_{n-1}^{-1/3})}
\E_x\exp\bigg\{\theta\int_0^{t_{n-1}}\ol{V}(B_s)ds\bigg\}dx
\end{aligned}
$$
when $t_n\le t\le t_{n+1}$ for large $n$. In view of (\ref{lb-10}) 
and (\ref{lb-17}), this leads to  (\ref{lb-13}). \qed

\section{Upper bounds}\label{ub}

In this section we install the upper bounds requested by
Theorem \ref{intro-4}, Theorem \ref{intro-9}, and Theorem \ref{intro-11}.
Through this section $0<\theta<1/16$.
Recall that $k=[(8\theta)^{-1}]$ and that for any open set $D\subset\R^3$
and the function $\zeta(\cdot)$ on $D$,
$\lambda_\zeta(D)$ is defined by the variation given in (\ref{FK-0'}).
For each $R>0$, write $Q_R=(-R, R)^d$.

\subsection{Asymptotics for the principal eigenvalues}\label{ub-1}

By $0<\theta<1/16$ we have that $k=[(8\theta)^{-1}]\ge 2$. Write
\begin{align}\label{ub-2}
R_k(t)=\left\{\begin{array}{ll} t^{k\over k-2}l(t)^{2\over 3(k-2)}
\hskip.4in \hbox{when}\hskip.1in k\ge 3\\\\
t^3l(t)^{2/3}\hskip.75in \hbox{when}\hskip.1in k=2.\end{array}\right.
\end{align}

\begin{lemma}\label{ub-3} 
$$
\lim_{t\to\infty}t^{-{2\over k-1}}l(t)^{-{2\over 3(k-1)}}\lambda_{\theta\ol{V}}
\big(Q_{R_k(t)}\big)=0\hskip.2in a.s.
$$
under the assumption
\begin{align}\label{ub-4}
\int_1^\infty{dt\over t\cdot l(t)}<\infty.
\end{align}
\end{lemma}

\proof We first consider the case $k\ge 3$. Let
$M>0$ be fixed but arbitrary.
Write
$$
r(t)=M\big(t l(t)^{1/3}\big)^{1\over (k-1)(k-2)},\hskip.2in
\epsilon(t)=\big(t^3 l(t)\big)^{-{k\over (k-1)(k-2)}}
$$
$$
\delta(t)=\epsilon(t)^{1/3}r(t)=M\big(t l(t)^{1/3}
\big)^{-{1\over k-2}}.
$$
Decompose $\ol{V}$ as follows:
$$
\begin{aligned}
\ol{V}(x)&=\int_{\R^3}{\alpha\big(\delta(t)^{-1}\vert y-x\vert\big)
\over\vert y-x\vert^2}\big[\omega(dy)-dy\big]
+\int_{\R^3}{1-\alpha\big(\delta(t)^{-1}\vert y-x\vert\big)
\over\vert y-x\vert^2}\big[\omega(dy)-dy\big].
\end{aligned}
$$
For the first term
$$
\begin{aligned}
\int_{\R^3} & {\alpha\big(\delta(t)^{-1}\vert y-x\vert\big)
\over\vert y-x\vert^2}\big[\omega(dy)-dy\big]
\le\int_{\R^3}{\alpha\big(\delta(t)^{-1}\vert y-x\vert\big)
\over\vert y-x\vert^2}\omega(dy)\\
&=\epsilon(t)^{-2/3}\int_{\R^3}{\alpha\big(r(t)^{-1}\vert y-
\epsilon(t)^{-1/3}x\vert\big)
\over\vert y-\epsilon(t)^{-1/3}x\vert^2}\omega(\epsilon(t)dy)
=\epsilon(t)^{-2/3}\xi_{r,\epsilon}\big(\epsilon(t)^{-1/3}x\big),
\end{aligned}
$$
where 
$$
\xi_{r,\epsilon}(x)=\xi_{r(t),\epsilon(t)}(x)
=\int_{\R^3}{\alpha\big(r(t)^{-1}\vert y-
x\vert\big)
\over\vert y-x\vert^2}\omega(\epsilon(t)dy).
$$

As for the second term
$$
\begin{aligned}
\int_{\R^3} & {1-\alpha\big(\delta(t)^{-1}\vert y-x\vert\big)
\over\vert y-x\vert^2}\big[\omega(dy)-dy\big]\\
&=a^2\delta(t)^{-2}\int_{\R^3}{1-\alpha\big(a^{-1}
\vert y-a\delta(t)^{-1}x\vert\big)
\over\vert y-a\delta(t)^{-1}x\vert^2}\big[\omega(a^{-3}\delta(t)^3dy)
-a^{-3}\delta(t)^3dy\big]\\
&=a^2\delta(t)^{-2}\ol{V}_{a, \tilde{\delta}(t)}\big(a\delta(t)^{-1}x\big)
\end{aligned}
$$
where $\tilde{\delta}(t)=a^{-3}\delta(t)^3$,
the random field $\ol{V}_{a,\epsilon}(\cdot)$ is defined in
(\ref{poisson-2'}) and the constant  $a>0$ will be specified later.

By triangle inequality and by the
substitution $g(x)\mapsto \epsilon(t)^{-1/2}g\Big(x\epsilon(t)^{-1/3}\Big)$,
\begin{align}\label{ub-5}
\lambda_{\theta\ol{V}}(Q_{R_k(t)})\le
\epsilon(t)^{-2/3}\lambda_{\theta\xi_{r,\epsilon}}
(Q_{\epsilon^{-1/3}(t)R_k(t)})+\theta a^2\delta^{-2}(t)
\sup_{x\in a\delta(t)^{-1} Q_{R_k(t)}}\vert
\ol{V}_{a, \tilde{\delta}(t)}(x)\vert.
\end{align}

By Proposition 1 in \cite{GK},  there is
a non-negative and continuous function $\Phi(x)$ on $\R^3$
whose support is contained in the 1-neighborhood of the grid $2r(t)\Z^3$,
such that
$$
\lambda_{\xi_{r,\epsilon}-\Phi^y}(Q_{\epsilon^{-1/3}(t)R_k(t)})
\le\max_{z\in 2r(t)\Z^3\cap  Q_{2\epsilon^{-1/3}(t)R_k(t)+2r(t)}}
\lambda_{\xi_{r,\epsilon}}(z+Q_{r(t)+1})
\hskip.2in y\in Q_{r(t)}
$$
where $\Phi^y(x)=\Phi(x+y)$.
In addition,  $\Phi(x)$ is periodic with period $2r(t)$:
$$
\Phi(x+2r(t)z)=\Phi(x);\hskip.2in x\in\R^3,\hskip.1in z\in\Z^3
$$
and there is a constant $K>0$ independent of $r(t)$ and $t$ such that
$$
\int_{Q_{r}}\Phi(x)dx\le {K\over r(t)}.
$$

By periodicity, therefore,
$$
\eta(x)\equiv{1\over \big(2r(t)\big)^3}
\int_{Q_{r}}\Phi^y(x)dy={1\over \big(2r(t)\big)^3}\int_{Q_{r}}\Phi(y)dy\le 
{K\over 8r(t)^4}.
$$
So we have
$$
\begin{aligned}
\lambda_{\xi_{r,\epsilon}}(Q_{\epsilon^{-1/3}(t)R_k(t)})
& \le {K\over 8r(t)^4}+\lambda_{\xi_{r,\epsilon}-\eta}(Q_{\epsilon^{-1/3}(t)R_k(t)})\\
&\le{K\over 8r(t)^4}+{1\over \big(2r(t)\big)^3}
\int_{Q_{r(t)}}
\lambda_{\xi_{r,\epsilon}-\Phi^y}(Q_{\epsilon^{-1/3}(t)R_k(t)})dy\\
&\le{K\over 8r(t)^4}+\max_{z\in 2r(t)\Z^3\cap  Q_{2\epsilon^{-1/3}(t)R_k(t)+2r(t)}}
\lambda_{\xi_{r,\epsilon}}(z+Q_{r(t)+1}),
\end{aligned}
$$
where the second step follows from Jensen inequality.

Summarizing the estimate since (\ref{ub-5}),
\begin{align}\label{ub-6}
t^{-{2\over k-1}}l(t)^{-{2\over 3(k-1)}} & \lambda_{\theta\ol{V}}
\big(Q_{R_k(t)}\big)\\
& \hskip-0.1in \le \theta a^2M^{-2}\sup_{x\in a\delta(t)^{-1} Q_{R_k(t)}}
\vert\ol{V}_{a, \tilde{\delta}(t)}(x)\vert+{K\over 8M^4}\nonumber\\
&+\big(tl(t)^{1/3}\big)^{2\over (k-1)(k-2)}
\max_{z\in 2r(t)\Z^3\cap  Q_{2\epsilon^{-1/3}(t)R_k(t)+2r(t)}}
\lambda_{\theta\xi_{r,\epsilon}}(z+Q_{r(t)+1}).\nonumber
\end{align}

Take $t_n=2^n$. By (\ref{poisson-10}),
\begin{align}\label{ub-7}
\lim_{n\to\infty}
\sup_{x\in a\delta(t_n)^{-1} Q_{R_k(t_n)}}
\vert\ol{V}_{a, \tilde{\delta}(t_n)}(x)\vert
=0\hskip.2in a.s.
\end{align}
when $a$ is sufficiently large. 

We now prove that
\begin{align}\label{ub-8}
\P\Big\{\max_{z\in 2r(t_n)\Z^3\cap  Q_{2\epsilon^{-1/3}(t_n)R_k(t_n)+2r(t_n)}}
\lambda_{\xi_{r(t_n),\epsilon(t_n)}}(z+Q_{r(t_n)+1})
=0\hskip.05in\hbox{eventually in $n$}
\Big\}=1.
\end{align}

Notice that
$$
\begin{aligned}
\P & \Big\{\max_{z\in 2r(t_n)\Z^3\cap  Q_{2\epsilon^{-1/3}(t_n)R_k(t_n)+2r(t_n)}}
\lambda_{\theta\xi_{r(t_n),\epsilon(t_n)}}(z+Q_{r(t_n)+1})\not =0\Big\}\\
&\le\#\Big\{2r(t_n)\Z^3\cap  Q_{2\epsilon^{-1/3}(t_n)R_k(t_n)+2r(t_n)}\Big\}
\P\Big\{\lambda_{\theta\xi_{r(t_n),\epsilon(t_n)}}(Q_{r(t_n)+1})\not =0\Big\}.
\end{aligned}
$$

Recall that truncation function
$\alpha (\cdot)$ is supported on $[0,3]$. 
For any $g\in {\cal F}_3(Q_{r(t_n)+1})$, 
$$
\begin{aligned}
\int_{Q_{r(t_n)+1}}\xi_{r(t_n),\epsilon(t_n)}(x) & g^2(x)  \, dx =\int_{\R^3}
\bigg[\int_{Q_{r_n(t)+1}}{\alpha(r(t_n)^{-1}\vert y-x\vert)\over\vert y-x\vert^2}
g^2(y)dy
\bigg]\omega\big(\epsilon(t_n)dx\big)\\
&=\int_{Q_{5r(t_n)}}\bigg[\int_{Q_{r(t_n)+1}}{\alpha(r(t_n)^{-1}\vert y-x\vert)
\over\vert y-x\vert^2}g^2(y)dy
\bigg]\omega\big(\epsilon(t_n)dx\big)\\
&\le\omega(Q_{5\delta(t_n)})\sup_{x\in \R^3}
\int_{Q_{r(t_n)+1}}{g^2(y)\over\vert y-x\vert^2}dy
\end{aligned}
$$
when $r(t_n)\ge 1$. Therefore,
\begin{align}\label{ub-8'}
&\lambda_{\theta\xi_{r(t_n),\epsilon(t_n)}}(Q_{r(t_n)+1})\\
&\le\sup_{g\in {\cal F}_3(Q_{r(t_n)+1})}\bigg\{
\omega(Q_{5\delta(t_n)})\theta\sup_{x\in \R^3}
\int_{Q_{r(t_n)+1}}{g^2(y)\over\vert y-x\vert^2}dy-{1\over 2}\int_{Q_{r_n(t)+1}}
\vert\nabla g(y)\vert^2dy\bigg\}\nonumber\\
&\le\sup_{g\in {\cal F}_3}\bigg\{
\omega(Q_{5\delta(t_n)})\theta\sup_{x\in \R^3}
\int_{\R^3}{g^2(y)\over\vert y-x\vert^2}dy-{1\over 2}\int_{\R^3}
\vert\nabla g(y)\vert^2dy\bigg\}\nonumber\\
&=\sup_{x\in\R^3}\sup_{g\in {\cal F}_3}\bigg\{
\omega(Q_{5\delta(t_n)})
\theta\int_{\R^3}{g^2(y)\over\vert y-x\vert^2}dy-{1\over 2}\int_{\R^3}
\vert\nabla g(y)\vert^2dy\bigg\}\nonumber\\
&=\sup_{g\in {\cal F}_3}\bigg\{
\omega(Q_{5\delta(t_n)})\theta
\int_{\R^3}{g^2(y)\over\vert y\vert^2}dy-{1\over 2}\int_{\R^3}
\vert\nabla g(y)\vert^2dy\bigg\}, \nonumber
\end{align}
where the last step follows from shifting invariance.

Notice that $k\theta\le 8^{-1}$. By Lemma \ref{H-2} we obtain the bound
$$
\begin{aligned}
\P & \Big\{\max_{z\in 2r(t_n)\Z^3\cap  Q_{2\epsilon^{-1/3}(t_n)R_k(t_n)+2r(t_n)}}
\lambda_{\theta\xi_{a,\epsilon(t_n)}}(z+Q_{r(t_n)+1})\not =0\Big\}\\
&\le C t_n^{3(k+1)\over k-2}l(t_n)^{3\over k-2}\P\Big\{
\omega(Q_{5\delta(t_n)})\ge k+1\Big\}\le Cl(t_n)^{-1}.
\end{aligned}
$$
By (\ref{ub-4}),
$$
\sum_nl(t_n)^{-1}<\infty.
$$
Hence, (\ref{ub-8}) follows from  Borel-Cantelli lemma.

Since the second term in (\ref{ub-6}) can be arbitrarily small by
making $M$ sufficiently large, by (\ref{ub-7}) and (\ref{ub-8}),
$$
\limsup_{n\to\infty}t_n^{-{2\over k-1}}l(t_n)^{-{2\over 3(k-1)}}
\lambda_{\theta\ol{V}}\big(Q_{R_k(t_n)}\big)\le 0\hskip.2in a.s.
$$
Notice that $\lambda_{\theta\ol{V}}\big(Q_{R_k(t)}\big)$ is
non-decreasing in $t$. We have completed the proof
in the case $k\ge 3$.

The case $k=2$ follows from the same argument with
$$
r(t)=M\big(tl(t)^{1/3}\big)^{1/2},
\hskip.2in \epsilon(t)=\big(t^3l(t)\big)^{-2},
$$
$$
\delta(t)=\epsilon(t)^{1/3}r(t)=M\big(tl(t)^{1/3}\big)^{-3/2}.
$$
\qed

Write
\begin{align}\label{ub-9}
S_k(t)=\left\{\begin{array}{ll} t^{k\over k-2}l(t)^{-{2\over 3(k-2)}}
\hskip.4in \hbox{when}\hskip.1in k\ge 3\\\\
t^3l(t)^{-2/3}\hskip.75in \hbox{when}\hskip.1in k=2.\end{array}\right.
\end{align}

\begin{lemma}\label{ub-10} 
$$
\liminf_{t\to\infty}t^{-{2\over k-1}}l(t)^{{2\over 3(k-1)}}\lambda_{\theta\ol{V}}
\big(Q_{S_k(t)}\big)=0\hskip.2in a.s.
$$
under the assumption that there is $c_0>0$ such that
\begin{align}\label{ub-12}
\int_1^\infty {1\over t}\exp\{-cl(t)\} \, dt \ 
\left\{\begin{array}{ll} =\infty 
\hskip.4in \hbox{when}\hskip.1in c<c_0\\\\
<\infty\hskip.4in \hbox{when}\hskip.1in c>c_0.\end{array}\right.
\end{align}
\end{lemma}

\proof  We first consider the case $k\ge 3$. Let
$u>0$ and $M>0$ be fixed but arbitrary.
Write
$$
r(t)=M\big(t l(t)^{-1/3}\big)^{1\over (k-1)(k-2)},\hskip.2in
\epsilon(t)=u^{-3}\big(t^3 l(t)^{-1}\big)^{-{k\over (k-1)(k-2)}}
$$
$$
\delta(t)=\epsilon(t)^{1/3}r(t)
=\Big({M\over u}\Big)\big(tl(t)^{-1/3}\big)^{-{1\over k-2}}.
$$

Similar to (\ref{ub-6}),
\begin{align}\label{ub-13}
t^{-{2\over k-1}}l(t)^{{2\over 3(k-1)}} & \lambda_{\theta\ol{V}}
\big(Q_{S_k(t)}\big)\\
& \hskip-0.1in \le \theta a^2\Big({u\over M}\Big)^2\sup_{x\in a\delta(t)^{-1} Q_{R_k(t)}}
\vert\ol{V}_{a, \tilde{\delta}(t)}(x)\vert+{Ku^2\over 8M^4}\nonumber\\
&+u^2\big(tl(t)^{-1/3}\big)^{2\over (k-1)(k-2)}
\max_{z\in 2r(t)\Z^3\cap  Q_{2\epsilon^{-1/3}(t)S_k(t)+2r(t)}}
\lambda_{\theta\xi_{r,\epsilon}}(z+Q_{r(t)+1}), \nonumber
\end{align} 
where $\tilde{\delta}=\tilde{\delta}(t)=a^{-3}\delta(t)^3$,
the random field $\ol{V}_{a,\epsilon}(x)$ is defined in (\ref{poisson-2'}),
and
$$
\xi_{r,\epsilon}(x)=\xi_{r(t),\epsilon(t)}(x)
=\int_{\R^3}{\alpha\big(r(t)^{-1}\vert y-
x\vert\big) \over\vert y-x\vert^2}\omega(\epsilon(t)dy).
$$

Same as (\ref{ub-8'}),
$$
\lambda_{\theta\xi_{r,\epsilon}}(z+Q_{r(t)+1})
\le\sup_{g\in {\cal F}_3}\bigg\{
\omega\big(\epsilon(t)^{1/3}(z+Q_{5r(t)})\big)\theta
\int_{\R^3}{g^2(y)\over\vert y\vert^2}dy-{1\over 2}\int_{\R^3}
\vert\nabla g(y)\vert^2dy\bigg\}
$$
for each $z\in 2r(t)\Z^d\cap  Q_{2\epsilon^{-1/3}(t)S_k(t)+2r(t)}$. Thus,
$$
\begin{aligned}
&\max_{z\in 2r(t)\Z^3\cap  Q_{2\epsilon^{-1/3}(t)S_k(t)+2r(t)}}
\lambda_{\theta\xi_{r,\epsilon}}(z+Q_{r(t)+1})\\
&\le\sup_{g\in {\cal F}_3}\bigg\{
\theta\Big(\max_{z\in 2r(t)\Z^3\cap  Q_{2\epsilon^{-1/3}(t)S_k(t)+2r(t)}}\omega
\big(\epsilon(t)^{1/3}(z+Q_{5r(t)})\big)\Big)
\int_{\R^3}{g^2(y)\over\vert y\vert^2}dy\\
&\hskip.6in -{1\over 2}\int_{\R^3}
\vert\nabla g(y)\vert^2dy\bigg\}\\
&=\sup_{g\in {\cal F}_3}\bigg\{
\theta\Big(\max_{z\in 2\delta(t)\Z^3\cap  Q_{2S_k(t)+2\delta(t)}}\omega
(z+Q_{5\delta(t)})\Big)
\int_{\R^3}{g^2(y)\over\vert y\vert^2}dy-{1\over 2}\int_{\R^3}
\vert\nabla g(y)\vert^2dy\bigg\}.
\end{aligned}
$$
By Lemma \ref{H-2}, therefore,
\begin{align}\label{ub-13'}
&\Big\{\max_{z\in 2r(t)\Z^3\cap  Q_{2\epsilon^{-1/3}(t)S_k(t)+2r(t)}}
\lambda_{\xi_{r,\epsilon}}(z+Q_{r(t)+1})=0\Big\}\\
& \quad \supset\Big\{\max_{z\in 2\delta(t)\Z^3\cap  Q_{2S_k(t)+2\delta(t)}}\omega
(z+Q_{5\delta(t)})\le k\Big\}.\nonumber
\end{align}

Unfortunately, the random variables
$$
\omega
(z+Q_{5\delta(t)})\hskip.2in 
z\in 2\delta(t)\Z^3\cap  Q_{2S_k(t)+2\delta(t)}
$$
are not independent. 
So we apply Slepian-type domination (Lemma \ref{lem:00}):
$$
\begin{aligned}
\P\Big\{\max_{z\in 2\delta(t)\Z^3\cap  Q_{2S_k(t)+2\delta(t)}}\omega
(z+Q_{5\delta(t)}) & \le  k\Big\}\\
&\ge \bigg(\P\Big\{\omega(Q_{5\delta(t)})\le k\Big\}
\bigg)^{\#\{2\delta(t)\Z^3\cap  Q_{2S_k(t)+2\delta (t)}\}}.\nonumber
\end{aligned}
$$

It is straightforward to check that
$$
\P\Big\{\omega(Q_{5\delta(t)})\ge k+1\Big\}
\sim {(10u^{-1}M)^{3(k+1)}\over (k+1)!}\Big(t^3l(t)^{-1}\Big)^{-{k+1\over k-2}}
\hskip.2in (t\to\infty)
$$
and that
$$
\#\{2\delta(t)\Z^3\cap  Q_{2S_k(t)+2\delta (t)}\}
\sim \Big({u\over M}\Big)^3t^{3(k+1)\over k-2}l(t)^{-{3\over k-2}}
\hskip.2in (t\to\infty).
$$
Hence, there is a constant $C_k$ independent of $u$ and $M$ such that
\begin{align}\label{ub-14}
\P\Big\{\max_{z\in 2\delta(t)\Z^3\cap  Q_{2S_k(t)+2\delta(t)}}\omega
(z+Q_{5\delta(t)})\le k\Big\}
\ge\exp\Big\{-C_k\Big({M\over u}\Big)^{3k}l(t)\Big\}
\end{align}
for large $t$. In connection to (\ref{ub-13}), our strategy is to make 
$u^2/M^4$, $M/u$ sufficiently small, and to make $u$ and $M$ 
sufficiently large.

Fix a constant $\tilde{c}$ satisfying
$$
{k-1\over 3k}c_0<\tilde{c}<c_0.
$$
Define $\{t_n\}$ as following:
$$
t_1=1,\hskip.1in t_{n+1}=t_n\exp\Big\{\tilde{c} 
l(t_n)\Big\}\hskip.2in n=1,2,\cdots.
$$
By (\ref{poisson-10}),
\begin{align}\label{ub-15}
\lim_{n\to\infty}
\sup_{x\in Q_{\epsilon^{-1/3}(t_n)S_k(t_n)}}\vert\ol{V}_{a,\tilde{\delta}(t_n)}(x)\vert
=0\hskip.2in a.s.
\end{align}
for sufficiently large $a$.

We now prove that
\begin{align}\label{ub-16}
\P\bigg\{\max_{z\in 2\delta(t)\Z^3\cap  Q_{2S_k(t)+2\delta(t)}}\omega
(z+Q_{5\delta(t)})\le k
\hskip.05in \hbox{i.o.}\bigg\}
=1
\end{align}

Write
$$
H_n=\max_{z\in 2\delta(t)\Z^3\cap  Q_{2S_k(t)+2\delta(t)}}\omega
(z+Q_{5\delta(t)}),
$$
$$
A_n=Q_{2S_k(t_{n+1})+2\delta(t_{n+1})}\setminus
Q_{2S_k(t_n)+b\delta(t_n)},
$$
$$
Z_n=\max_{z\in 2\delta(t_n)\Z^3\cap A_n}\omega
(z+Q_{5\delta(t_{n+1})}),
$$
$$
\widetilde{Z}_n=\max_{z\in 2\delta(t_n)\Z^3\cap  Q_{2S_k(t_n)+b\delta(t_n)}}
\omega(z+Q_{5\delta(t_{n+1})}),
$$
where $b>0$ is a constant which is large enough to make sure that the
random variables $Z_1, Z_2,\cdots$ are independent.

We have that $H_{n+1}=\max\{Z_n, \widetilde{Z}_n\}$.
Notice that
$$
\begin{aligned}
\P\{\widetilde{Z}_n\ge k+1\}
&\le \#\{2\delta(t_n)\Z^3\cap  Q_{2S_k(t_n)+b\delta(t_n)}\}
\P\Big\{\omega(Q_{5\delta(t_{n+1})}\big)\ge k+1\Big\}\\
&\le Ct_n^{3(k+1)\over k-1}l(t_n)^{-{3\over k-2}}
t_{n+1}^{-{3(k+1)\over k-1}}l(t_{n+1})^{k+1\over k-2}\\
&=Cl(t_n)^{-{3\over k-2}}l(t_{n+1})^{k+1\over k-2}
\exp\Big\{-{3\tilde{c}(k+1)\over k-2}l(t_n)\Big\}.
\end{aligned}
$$
Since $l(t)$ is slow-varying,
$$
l(t_{n+1})=l\Big(t_n\exp\{\tilde{c}l(t_n)\}\Big)
\le l(t_n)\exp\Big\{o\big(l(t_n)\Big\}
=\exp\Big\{o\big(l(t_n)\Big\}
$$
for large $n$. Therefore, we obtain the bound
$$
\P\{\widetilde{Z}_n\ge k+1\}\le C\exp\Big\{-{3k\tilde{c}\over k-2}l(t_n)\Big\}
\hskip.2in (n\to\infty).
$$

For any $c>c_0$, on the other hand,
$$
\begin{aligned}
\infty &>\int_1^\infty {1\over t}\exp\big\{-cl(t)\big\}dt
=\sum_{n=1}^\infty\int_{t_n}^{t_{n+1}}{1\over t}\exp\big\{-cl(t)\big\}dt\\
&\ge \sum_{n=1}^\infty{t_{n+1}-t_n\over t_{n+1}}\exp\big\{-cl(t_{n+1})\big\}
\ge \delta\sum_{n=1}^\infty\exp\big\{-cl(t_{n+1})\big\}
\end{aligned}
$$
So we have that
$$
\sum_n\P\{\widetilde{Z}_n\ge k+1\}<\infty.
$$
By Borel-Cantelli lemma
\begin{align}\label{ub-16'}
\P\{\widetilde{Z}_n\le k\hskip.05in\hbox{eventually in $n$}\}=1.
\end{align}

By (\ref{ub-14}),
$$
\P\{Z_n\le k\}\ge\P\{H_{n+1}\le k\}\ge
\exp\Big\{-C_k\Big({M\over u}\Big)^{3k}l(t_{n+1})\Big\}.
$$

Pick $c_1$ satisfying $\tilde{c}<c_1<c_0$ and make $M/u$ so small that
$$
C_k\Big({r\over u}\Big)^{3k}<c_1-\tilde{c}
$$
We have
$$
\begin{aligned}
\infty &=\int_1^\infty {1\over t}\exp\big\{-c_1l(t)\big\}dt
=\sum_{n=1}^\infty\int_{t_n}^{t_{n+1}}{1\over t}\exp\big\{-c_1l(t)\big\}dt\\
&\le \sum_{n=1}^\infty{t_{n+1}-t_n\over t_{n}}\exp\big\{-c_1l(t_{n})\big\}
\le\sum_{n=1}^\infty\exp\big\{-(c_1 -\tilde{c})l(t_{n})\big\}.
\end{aligned}
$$
Consequently,
$$
\sum_n\P\{Z_n\le k\}=\infty.
$$
Applying Borel-Cantelli lemma to the independent sequence $\{Z_n\}$
we have
$$
\P\{Z_n\le k\hskip.05in\hbox{i.o.}\}=1.
$$
This, together with (\ref{ub-16'}), leads to
(\ref{ub-16}).
By (\ref{ub-13'}) and (\ref{ub-16}),
\begin{align}\label{ub-17}
\P\Big\{\max_{z\in 2r(t_n)\Z^3\cap  Q_{2\epsilon^{-1/3}(t_n)S_k(t_n)+2r(t_n)}}
\lambda_{\xi_{r(t_n),\epsilon(t_n)}}(z+Q_{r(t_n)+1})=0\hskip.1in\hbox{i.o.}\Big\}=1.
\end{align}

By (\ref{ub-13}), (\ref{ub-15}), (\ref{ub-17}), and by the fact that
$u^2/M^4$ can be arbitrarily small,
$$
\liminf_{n\to\infty}t_n^{-{2\over k-1}}l(t_n)^{{2\over 3(k-1)}}\lambda_{\theta\ol{V}}
\big(Q_{S_k(t_n)}\big)\le 0\hskip.2in a.s.
$$
By the fact that the principal eigenvalue $\lambda_{\theta\ol{V}}(Q_R)$
increases in $R$, we have completed the proof in the case $k\ge 3$.

The case $k=2$ follows from the same argument with
$$
r(t)=M\big(tl(t)^{-1/3}\big)^{1/2},
\hskip.2in \epsilon(t)=u^{-3}\big(t^3l(t)^{-1}\big)^{-2},
$$
$$
\delta(t)=\epsilon(t)^{1/3}r(t)=\Big({M\over u}\Big)\big(tl(t)^{-1/3}
\big)^{-3/2}.
$$
\qed

\subsection{Upper bound for Theorem \ref{intro-4}}

We prove that when $\theta<16^{-1}$,
\begin{align}\label{ub-18}
\E_0\exp\bigg\{\theta\int_0^t\ol{V}(B_s)ds\bigg\}<\infty\hskip.2in a.s.
\end{align}
for any $t>0$. By H\"older inequality, 
we may assume that $\theta> {1\over 24}$.

Let $l(t)\ge 0$ be a slow-varying function satisfying (\ref{ub-4})
and recall the notation $R_2(t)=t^3l(t)^{2/3}$.
Consider the decomposition
$$
\begin{aligned}
\E_0 & \exp\bigg\{\theta\int_0^t\ol{V}(B_s)ds\bigg\}
=\E_0\Bigg[\exp\bigg\{\theta\int_0^t\ol{V}(B_s)ds\bigg\};\hskip.05in
\tau_{Q_{R_2(t)}}\ge 2t\Bigg]\\
&+\sum_{n=1}^\infty\E_0\Bigg[\exp\bigg\{\theta\int_0^t\ol{V}(B_s)ds\bigg\};
\hskip.05in
\tau_{Q_{R_2(2^{n-1}t)}}<2t\le\tau_{Q_{R_2(2^{n}t)}} \Bigg].
\end{aligned}
$$

Pick $p>1$ with $p\theta<16^{-1}$ and write $q=p(p-1)^{-1}$. By H\"older
inequality,
$$
\begin{aligned}
\E_0 & \Bigg[\exp\bigg\{\theta\int_0^t\ol{V}(B_s)ds\bigg\};
\hskip.05in
\tau_{Q_{R_2(2^{n-1}t)}}<2t\le\tau_{Q_{R_2(2^{n}t)}} \Bigg]\\
&\le\Big(\P_0\big\{\tau_{Q_{R_2(2^{n-1}t)}}<2t\big\}\Big)^{1/q}
\Bigg\{\E_0\Bigg[\exp\bigg\{p\theta\int_0^t\ol{V}(B_s)ds\bigg\};
\hskip.05in\tau_{Q_{R_2(2^{n}t)}}\ge 2t \Bigg]\Bigg\}^{1/p}.
\end{aligned}
$$

Let $\delta>0$ be a small number and condition
on the event $\big\{\omega\big(B(0,\delta)\big)=0\big\}$. 
Applying Lemma \ref{FK-1}, 
\begin{align}\label{ub-19}
\E_0 & \exp\bigg\{\theta\int_0^t\ol{V}(B_s)ds\bigg\}
\le X(\theta t)+
Y_0(t)\exp\Big\{t\lambda_{\theta\ol{V}}(Q_{R_2(t)})\Big\}\\
&+\sum_{n=1}^\infty \Big(\P_0\big\{\tau_{Q_{R_2(2^{n-1}t)}}<2t\big\}\Big)^{1/q}
\bigg(X(p\theta t)+Y_n(t)
\exp\Big\{t\lambda_{p\theta\ol{V}}(Q_{R_2(2^{n}t)})\Big\}\bigg)^{1/p}, \nonumber
\end{align}
where
$$
X(t)=\exp\Big\{ t\sup_{\vert x\vert\le \delta/2}
\vert\ol{V}_{{\delta\over 6}, 1}(x)\vert\Big\},
$$
$$
Y_0(t)={48 R_2(t)^3\over\pi\delta^3} \, 
\E_0\exp\Big\{\sqrt{2\delta}\theta T_1\sup_{x\in Q_{R_2(t)}}
\vert\ol{V}_{{\delta\over 6}, 1}(x)\vert\Big\}
$$
and 
$$
Y_n(t)={48 R_2(2^nt)^3\over\pi\delta^3} \,  
\E_0\exp\Big\{\sqrt{2\delta}p\theta T_1\sup_{x\in Q_{R_2(2^nt)}}
\vert\ol{V}_{{\delta\over 6}, 1}(x)\vert\Big\}\hskip.2in n=1,2,\cdots.
$$

Using the classical fact that there is a constant $C>0$ such that
$$
\E_0\exp\{bT_1\}\le \exp\{Cb^2\}\hskip.2in \forall b>0
$$
we have
$$
\begin{aligned}
\E_0 & \exp\Big\{\sqrt{2\delta}p\theta T_1\sup_{x\in Q_{R_2(2^nt)}}
\vert\ol{V}_{{\delta\over 6}, 1}(x)\vert\Big\}
\le \exp\bigg\{2\delta C(p\theta)^2\Big(\sup_{x\in Q_{R_2(2^nt)}}
\vert\ol{V}_{{\delta\over 6}, 1}(x)\vert\Big)^2\bigg\}\\
&=\exp\Big\{o\Big((\log (2^nt)\big)^2\Big)\Big\}\hskip.2in a.s.
\hskip.2in (n\to\infty)
\end{aligned}
$$
where the last step follows from (\ref{poisson-9}). Consequently,
\begin{align}\label{ub-20}
Y_n(t)=\exp\Big\{o\big(n^2\big)\Big\}\hskip.2in a.s.
\hskip.2in (n\to\infty).
\end{align}

Recall the classic fact that
\begin{align}\label{ub-21}
\P_0 & \big\{\tau_{Q_{R_2(2^{n-1}t)}}<2t\big\}
=\P\Big\{\max_{s\le 2t}\vert B_s\vert_\infty\ge R_2(2^{n-1}t)\Big\}\\
&=\P\Big\{\max_{s\le 1}\vert B_s\vert_\infty\ge (2t)^{-1/2}R_2(2^{n-1}t)\Big\}
\le\exp\Big\{-C2^{6n}t^5l(2^{n-1}t)^{4/3}\Big\}\nonumber
\end{align}
for some constant $C>0$ independent of $n$ and $t$, where 
$\vert\cdot\vert_\infty$
is the max-norm in $\R^3$.

By Lemma \ref{ub-3} with $k=2$ and with $\theta$ being replaced by $p\theta$,
\begin{align}\label{ub-22}
\lambda_{p\theta\ol{V}}(Q_{R_2(2^{n}t)})
=o\Big(\big(2^nt\big)^2l\big(2^nt\big)^{2/3}\Big)\hskip.1in a.s.\hskip.2in
(n\to\infty).
\end{align}

Combining (\ref{ub-20}), (\ref{ub-21}),  (\ref{ub-22})
we conclude that the right hand side of (\ref{ub-19}) is almost
surely finite. Thus, we have established (\ref{ub-18}) conditioning
on the event $\Big\{\omega\big(B(0,\delta)\big)=0\Big\}$. Therefore,
$$
\P\bigg\{\E_0\exp\bigg\{\theta\int_0^t\ol{V}(B_s)ds\bigg\}<\infty\bigg\}
\ge \P\Big\{\omega\big(B(0,\delta)\big)=0\Big\}
=\exp\Big\{-{4\over 3}\pi\delta^3\Big\}.
$$
Since $\delta$ can be arbitrarily small, we have completed the proof. \qed

\subsection{Upper bound for Theorem \ref{intro-9}}

Consider the decomposition
$$
\ol{V}(x)=\ol{V}_{1,1}(x)
+\int_{\R^3}{\alpha(\vert y-x\vert)\over\vert y-x\vert^2}\omega(dy)
-\int_{\R^3}{\alpha(\vert y\vert)\over\vert y\vert^2}dy
\ge \ol{V}_{1,1}(x)-\int_{\R^3}{\alpha(\vert y\vert)\over\vert y\vert^2}dy
$$
where the notation $\ol{V}_{1,1}(x)$ comes from (\ref{poisson-2'}).
We have that
$$
\begin{aligned}
\E_0 & \exp\bigg\{\theta\int_0^t\ol{V}(B_s)ds\bigg\}
\ge \E_0\Bigg[\exp\bigg\{\theta\int_0^t\ol{V}(B_s)ds\bigg\};\hskip.05in
\tau_{B(0,t)}\ge t\Bigg]\\
&\ge\exp\bigg\{-t\bigg(\sup_{\vert x\vert\le t}\vert V_{1,1}(x)\vert
+\int_{\R^3}{\alpha(\vert y\vert)\over\vert y\vert^2}dy\bigg)\bigg\}
\P_0\Big\{\max_{s\le t}\vert B_s\vert\le t\Big\}.\nonumber
\end{aligned}
$$
By (\ref{poisson-9}) we have
\begin{align}\label{ub-23}
\liminf_{t\to\infty} \, (t\log t)^{-1}
\log\E_0\exp\bigg\{\theta\int_0^t\ol{V}(B_s)ds\bigg\}\ge 0\hskip.2in a.s.
\end{align}

To complete the proof of Theorem \ref{intro-9}, therefore, all we need to show
is that under the assumption (\ref{ub-4}),
\begin{align}\label{ub-24}
\limsup_{t\to\infty}t^{-{k+1\over k-1}}l(t)^{-{2\over 3(k-1)}}
\log\E_0\exp\bigg\{\theta\int_0^t\ol{V}(B_s)ds\bigg\}\le 0\hskip.2in a.s.
\end{align}
conditioning on the event $\Big\{\omega\big(B(0,\delta)\big)=0\Big\}$.

In the case $1/24 <\theta< 1/16$, the bound
(\ref{ub-19}) holds when conditioned on ${\Big\{\omega\big(B(0,\delta)\big)=0\Big\}}$.

By Lemma \ref{ub-3} with $k=2$,
$$
\lambda_{\theta\ol{V}}(Q_{R_2(t)})
=o\Big(t^2l(t)^{2/3}\Big)\hskip.1in a.s.\hskip.2in
(t\to\infty).
$$
The bound (\ref{ub-20}) can be replaced by
$$
Y_n(t)=\exp\Big\{o\Big(\big(\log(2^nt)\big)^2\Big)\Big\}\hskip.2in a.s.
\hskip.2in n=0,1,\cdots.
$$
Combining these with the bound given in  (\ref{ub-21}), 
(\ref{ub-22}), we have (\ref{ub-24}).

Now we consider the case $\displaystyle 0<\theta\le 1/24$ (so $k\ge 3$).
The main reason we treat this setting
separately is for it includes the critical cases when $\theta =(8k)^{-1}$,
which need some special care. Similar to (\ref{ub-19}), for any
conjugate $p, q>1$,
\begin{align}\label{ub-25}
\E_0 & \exp\bigg\{\theta\int_0^t\ol{V}(B_s)ds\bigg\}
\le X(\theta t)+
Y_0(t)\exp\Big\{t\lambda_{\theta\ol{V}}(Q_{R_k(t)})\Big\}\\
&+\sum_{n=1}^\infty \Big(\P_0\big\{\tau_{Q_{R_k(2^{n-1}t)}}<2t\big\}\Big)^{1/q}
\bigg(X(p\theta t)+Y_n(t)
\exp\Big\{t\lambda_{p\theta\ol{V}}(Q_{R_k(2^{n}t)})\Big\}\bigg)^{1/p}, \nonumber
\end{align}
where
\begin{align*}
X(t) & =\exp\Big\{ t\sup_{\vert x\vert\le \delta/2}
\vert\ol{V}_{{\delta\over 6}, 1}(x)\vert\Big\}, \\
Y_0(t) &= {48 R_k(t)^3\over\pi\delta^3}
\E_0\exp\Big\{\sqrt{2\delta}\theta T_1\sup_{x\in Q_{R_k(t)}}
\vert\ol{V}_{{\delta\over 6}, 1}(x)\vert\Big\}, \\
Y_n(t) &= {48 R_k(2^nt)^3\over\pi\delta^3}
\E_0\exp\Big\{\sqrt{2\delta}p\theta T_1\sup_{x\in Q_{R_k(2^nt)}}
\vert\ol{V}_{{\delta\over 6}, 1}(x)\vert\Big\}\hskip.2in n=1,2,\cdots.
\end{align*}
Similarly to (\ref{ub-20}) and (\ref{ub-21}) we get
$$
Y_n(t)=\exp\Big\{o\Big(\big(\log (2^nt)\big)^2\Big)\Big\}\hskip.2in a.s.
\hskip.2in n=0, 1,\cdots
$$
and
$$
\P_0\big\{\tau_{Q_{R_k(2^{n-1}t)}}<2t\big\}
\le\exp\Big\{-C2^n (2^nt)^{k+2\over k-2}l(2^{n-1}t)^{4\over 3(k-2)}\Big\}.
$$

Due to the possibility that $\theta =(8k)^{-1}$, we can only
make $p\theta <\big(8(k-1)\big)^{-1}$. So we
may make $(8k)^{-1}<p\theta <\big(8(k-1)\big)^{-1}$. By
the monotonicity of $\lambda_{p\theta \ol{V}}(D)$ in $D$,
$$
\lambda_{p\theta\ol{V}}(Q_{R_k(2^nt)})\le \lambda_{p\theta\ol{V}}(Q_{R_{k-1}(2^nt)})
=o\bigg(\Big(2^nt\Big)^{2\over k-2}l\Big(2^n t\Big)^{2\over 3(k-2)}
\bigg)\hskip.2in a.s.
$$
where the second step follows
from Lemma \ref{ub-3} with $k$ being replaced by $k-1$.

Summarizing the bounds we obtained, the infinite series on 
the right hand side of (\ref{ub-25}) is asymptotically (as $t\to\infty$)
and almost surely
bounded by 
$$
C\sum_{n=1}^\infty \exp\Big\{-C^{-1}2^n\Big\}.
$$

We now obtain desired (\ref{ub-24}) applying (\ref{ub-25}) and and the fact  that
$$
X(\theta t)+
Y_0(t)\exp\Big\{t\lambda_{\theta\ol{V}}(Q_{R_k(t)})
=\exp\Big\{o\Big(t^{k+1\over k-1}l(t)^{2\over 3(k-1)}\Big)\Big\}\hskip.2in a.s.
\hskip.2in (t\to\infty).
$$
(see Lemma \ref{ub-3}). The proof is complete.
\qed

\subsection{Upper bound for Theorem \ref{intro-11}}

In view of (\ref{ub-23}),  we only need to show
\begin{align}\label{ub-26}
\liminf_{t\to\infty}t^{-{k+1\over k-1}}l(t)^{{2\over 3(k-1)}}
\log\E_0\exp\bigg\{\theta\int_0^t\ol{V}(B_s)ds\bigg\}\le 0\hskip.2in a.s.
\end{align}
conditioning on the event $\Big\{\omega\big(B(0,\delta)\big)=0\Big\}$.

We prove (\ref{ub-26}) under the extra assumption that
$$
\int_1^\infty {1\over t}\exp\big\{-cl(t)\big\}dt<\infty
$$
for some large constant $c>0$, for otherwise we may consider
$\tilde{l}(t)=\log\log t+l(t)$ instead of $l(t)$. Therefore,
(\ref{ub-12}) can be assumed here.

Let $S_k(t)$ be given as in Lemma \ref{ub-10}. We have that
\begin{align}\label{ub-27}
\E_0 & \exp\bigg\{\theta\int_0^t\ol{V}(B_s)ds\bigg\}\\
&\le\E_0\Bigg[\exp\bigg\{\theta\int_0^t\ol{V}(B_s)ds\bigg\};\hskip.05in 
\tau_{Q_{S_k(t)}}\ge 2t\Bigg]\nonumber\\
& \quad +\Big(\P_0\big\{\tau_{Q_{S_k(t)}}< 2t\big\}\Big)^{1/q}
\Bigg(\E_0\exp\bigg\{p\theta\int_0^t\ol{V}(B_s)ds\bigg\}\Bigg)^{1/p}\nonumber
\end{align}
where $p, q>1$ are conjugate numbers.

In the case $\displaystyle {1/24}<\theta <{1/16}$ ($k=2$),
we can make $p$ close to 1 so $\displaystyle p\theta<{1/ 16}$.
By the upper bound in Theorem \ref{intro-4} (with $\theta$ being
replaced by $p\theta$ and $l(t)=(\log t)^2$)
$$
\E_0\exp\bigg\{p\theta\int_0^t\ol{V}(B_s)ds\bigg\}
=\exp\bigg\{o\Big(t^3(\log t)^{4/3}\Big)\bigg\}\hskip.2in a.s.\hskip.2in
(t\to\infty).
$$

By the bound for Gaussian tail,
$$
\P_0\big\{\tau_{Q_{S_2(t)}}< 2t\big\}
\le \exp\Big\{-Ct^{-1}S_2(t)^2\Big\}=\exp\Big\{-Ct^5l(t)^{-4/3}\Big\}.
$$

Hence, the second term on the right hand side of (\ref{ub-27}) 
is negligible when 
$\displaystyle {1/24}<\theta <{1/16}$.

We now show that the same thing happens in the case when
$\displaystyle 0<\theta \le{1/24}$ ($k\ge 3$). 
In this case we can pick $p>1$ such that
$(8k)^{-1}<p\theta<\big(8(k-1)\big)^{-1}$. By Theorem \ref{intro-9}
(with $l(t)=(\log t)^2$ and $k$ being replaced by $k-1$),
$$
\E_0\exp\bigg\{p\theta\int_0^t\ol{V}(B_s)ds\bigg\}
=\exp\bigg\{o\Big(t^{k\over k-2}(\log t)^{4\over 3(k-2)}\Big)
\bigg\}\hskip.2in a.s.\hskip.2in
(t\to\infty).
$$
So our assertion follows from the Gaussian tail estimate
$$
\P_0\big\{\tau_{Q_{S_2(t)}}< 2t\big\}
\le \exp\Big\{-Ct^{-1}S_k(t)^2\Big\}=\exp\Big\{-Ct^{2k-1\over k-2}
l(t)^{-{4\over 3(k-2)}}\Big\}.
$$

Therefore, the problem (in both $k=2$ and $k\ge 3$)
has been reduced to the proof of
\begin{align}\label{ub-28}
\liminf_{t\to\infty}t^{-{k+1\over k-1}}l(t)^{{2\over 3(k-1)}}
\log\E_0\Bigg[\exp\bigg\{\theta\int_0^t\ol{V}(B_s)ds\bigg\};
\hskip.05in \tau_{Q_{S_k(t)}}\ge 2t\Bigg]\le 0\hskip.2in a.s.
\end{align}
conditioning on the event $\Big\{\omega\big(B(0,\delta)\big)=0\Big\}$.

By Lemma \ref{FK-1},
$$
\begin{aligned}
\E_0 & \Bigg[\exp\bigg\{\theta\int_0^t\ol{V}(B_s)ds\bigg\};
\hskip.05in \tau_{Q_{S_k(t)}}\ge 2t\Bigg]
\le\exp\Big\{\theta t\sup_{\vert x\vert\le\delta/2}\vert\ol{V}_{{\delta\over 6}, 1}
(x)\vert\Big\}\\
&+{6\vert Q_{S_k(t)}\vert\over\pi\delta^3}
\E_0\exp\Big\{\sqrt{2\delta}T_1\theta\sup_{x\in Q_{S_k(t)}}
\vert\ol{V}_{{\delta\over 6}, 1}
(x)\vert\Big\}\exp\Big\{t\lambda_{\theta\ol{V}}(Q_{S_k(t)})\Big\}\\
&=\exp\{O(t)\}+\exp\Big\{\Big(o\big(\log S_k(t)\big)^2\Big)\Big\}
\exp\Big\{t\lambda_{\theta\ol{V}}(Q_{S_k(t)})\Big\}\hskip.2in a.s.
\hskip.2in (t\to\infty)
\end{aligned}
$$
where the last step follows from (\ref{poisson-9}).

The required (\ref{ub-28}) follows from Lemma \ref{ub-10}. \qed

\section{Hardy inequality}\label{H}

Recall the definition of ${\cal F}_d(D)$
from (\ref{FK-0}).  The family ${\cal F}_3$ is defined as
$$
{\cal F}_3={\cal F}_3(\R^3)=\bigg\{g\in W^{1,2}(\R^3);\hskip.1in
\int_{\R^3}g^2(x)dx=1\bigg\}.
$$

The essential reason behind the main theorems in this paper is the
Hardy's inequality. Searching in literature, we have found large
amount of follow-up publication (i.e., \cite{HPL} and \cite{OK}) on
this subject, except Hardy's original paper. For reader's convenience,
we state  Hardy's inequality for $d=3$ in the following lemma and
provide a short proof.

\begin{lemma}\label{H-0'} For any $f_\epsilon\in W^{1,2}(\R^3)$,
\begin{align}\label{H-0}
\int_{\R^3}{f^2(x)\over\vert x\vert^2}dx\le 4\int_{\R^3}\vert\nabla f(x)\vert^2
dx.
\end{align}
Further, the number 
4 is the best constant in the sense that for any $\epsilon>0$ one
can find a function $f_\epsilon\in W^{1,2}(\R^3)$ with compact support such that
\begin{align}\label{H-1}
\int_{\R^3}{f_\epsilon^2(x)\over\vert x\vert^2}dx> (4-\epsilon)
\int_{\R^3}\vert\nabla f_\epsilon(x)\vert^2
dx.
\end{align}
\end{lemma}

\proof Write $x=(x_1,x_2,x_3)$. Using integration by parts
$$
\int_{\R^3}{f^2(x)\over\vert x\vert^2}dx
=\int_{\R^3}x_j\Big[{2x_i\over\vert x\vert^4}f^2(x)-{2\over\vert x\vert^2}
f(x){\partial f\over\partial x_j}\Big]dx\hskip.2in j=1,2,3.
$$
Summing over $j$ on the both sides
$$
3\int_{\R^3}{f^2(x)\over\vert x\vert^2}dx=
2\int_{\R^3}\Big[{f^2(x)\over\vert x\vert^2}-
{\nabla f\cdot x\over\vert x\vert^2}f(x)\Big]dx.
$$
Thus,
$$
\int_{\R^3}{f^2(x)\over\vert x\vert^2}dx
=-2\int_{\R^3}{\nabla f\cdot x\over\vert x\vert}{f(x)\over\vert x\vert}dx
\le 2\bigg(\int_{\R^3}{\vert \nabla f\cdot x\vert^2\over\vert x\vert^2}
dx\bigg)^{1/2}\bigg(\int_{\R^3}{f^2(x)\over\vert x\vert^2}dx\bigg)^{1/2}.
$$
Therefore,
$$
\int_{\R^3}{f^2(x)\over\vert x\vert^2}dx\le 4
\int_{\R^3}{\vert \nabla f\cdot x\vert^2\over\vert x\vert^2}dx
\le 4\int_{\R^3}\vert \nabla f(x)\vert^2dx.
$$

To establish (\ref{H-1}), for each large $M>0$, we define 
$g_M\in W^{1,2}(\R^3)$
as following:
$$
\begin{aligned}
g_M(x)=\left\{\begin{array}{ll}M^{1/2}\hskip.3in 0\le\vert x\vert\le M^{-1}\\\\
\vert x\vert^{-1/2}\hskip.3in M^{-1}<\vert x\vert\le M\\\\
\displaystyle{2M-\vert x\vert\over M^{3/2}}\hskip.3in M<\vert x\vert\le 2M\\\\
0\hskip.9in\vert x\vert> 2M.\end{array}\right.
\end{aligned}
$$
It is straightforward to exam that $g_M$ is locally supported and
$$
\int_{\R^3}{g_M^2(x)\over\vert x\vert^2}dx=
\bigg\{4-28\Big({7\over 3}+{1\over 2}\log M\Big)^{-1}\bigg\}
\int_{\R^3}\vert\nabla g_M(x)\vert^2dx.
$$
For each $\epsilon>0$, take $M>0$ sufficiently large so
$$
28\Big({7\over 3}+{1\over 2}\log M\Big)^{-1}<\epsilon
$$
and let $f_\epsilon(x)=g_M(x)$.\qed

What has been frequently used in this paper is the following version
of Hardy's inequality.

\begin{lemma}\label{H-2} For any $\theta>0$,
\begin{align}\label{H-3}
\sup_{g\in {\cal F}_3}\bigg\{\theta\int_{\R^3}{g^2(x)\over\vert x\vert^2}dx
-{1\over 2}\int_{\R^3}\vert\nabla g(x)\vert^2dx\bigg\}
= \begin{cases}
 0   & \text{if } \theta\le 1/8, \medskip \\ 
 \infty   & \text{if } \theta > 1/8. 
\end{cases}
\end{align}
\end{lemma}

\proof By Hardy's inequality, the left hand side of (\ref{H-3}) 
is non-positive when $\theta<1/8$. On the other hand, it is no less than
$$
-{1\over 2}\inf_{g\in {\cal F}_3}\int_{\R^3}\vert\nabla g(x)\vert^2dx
$$
which is equal to zero. Thus, for $\theta\le 1/8$,
$$
\sup_{g\in {\cal F}_3}\bigg\{\theta\int_{\R^3}{g^2(x)\over\vert x\vert^2}dx
-{1\over 2}\int_{\R^3}\vert\nabla g(x)\vert^2dx\bigg\}=0.
$$

Assume $\theta> 1/8$. By the optionality of Hardy's inequality described
in (\ref{H-1}),
$$
H(\theta)\equiv\sup_{g\in {\cal F}_3}
\bigg\{\theta\int_{\R^3}{g^2(x)\over\vert x\vert^2}dx
-{1\over 2}\int_{\R^3}\vert\nabla g(x)\vert^2dx\bigg\}>0.
$$
Given $a>0$, the substitution $g(x)=a^{3/2}f(ax)$ leads to
$H(\theta)=a^2H(\theta)$. So $M(\theta)=\infty$. \qed

\bigskip	

\section*{Acknowledgments}

{\small We wish to thank Witold Nazarewicz of Physics Department at University of Tennessee for helpful discussions and suggesting relevant literature in physics on the subject.}

\bibliographystyle{amsplain}

\end{document}